\documentclass{article}

\usepackage{amsmath, amscd, amsthm, amssymb, graphics} 

\newcommand{\liets}{{\mathfrak t}^*}
\newcommand{\Proj}{{\rm Proj}}
\newcommand{\Spec}{{\rm Spec}}
\newcommand{\diag}{{\rm diag}}

\input amssym.def
\input amssym.tex
\newcommand{\nc}{\newcommand}
\nc{\bla}{\phantom{bbbbb}}

\newcommand{\beq}{\begin{equation}}
\newcommand{\eeq}{\end{equation}}
\newcommand{\barr}{\begin{array}}
\newcommand{\earr}{\end{array}}
\newcommand{\beqar}{\begin{eqnarray}}
\newcommand{\eeqar}{\end{eqnarray}}
\newtheorem{theorem}{Theorem}[section]
\newtheorem{corollary}[theorem]{Corollary}
\newtheorem{lemma}[theorem]{Lemma}
\newtheorem{prop}[theorem]{Proposition}
\newtheorem{definition}[theorem]{Definition}
\newtheorem{remit}[theorem]{Remark}

\newtheorem{exit}[theorem]{Example}

\newenvironment{rem}{\begin{remit}\rm}{\end{remit}}
\newenvironment{ex}{\begin{exit}\rm}{\end{exit}}
\newenvironment{defn}{\begin{definition}\rm}{\end{definition}}

\newcommand{\RR}{{\mathbb R }}
\newcommand{\CC}{{\mathbb C }}
\nc{\FF}{ {\mathbb F} } 
\nc{\HH}{ {\mathbb H} } 
\newcommand{\ZZ}{{\mathbb Z }}
\newcommand{\PP}{ {\mathbb P } }
\newcommand{\QQ}{{\mathbb Q }}

\newcommand{\calf}{{\mbox{$\cal F$}}}

\newcommand{\cali}{{\mbox{$\cal I$}}}

\newcommand{\calm}{{\mbox{$\cal M$}}}

\newcommand{\calo}{{\mbox{$\cal O$}}}
\newcommand{\calp}{{\mbox{$\cal P$}}}

\newcommand{\calx}{{\mbox{$\cal X$}}}
\newcommand{\caly}{{\mbox{$\cal Y$}}}
\newcommand{\calz}{{\mbox{$\cal Z$}}}

\nc{\conv}{{\rm Conv}}
\nc{\umax}{{U_{\max}}}

\newcommand{\liek}{{\mathfrak k}}
\newcommand{\lieu}{{\mathfrak u}}

\newcommand{\lieks}{{\liek}^*}
\newcommand{\liet}{{\mathfrak t}}
\newcommand{\xg}{X/\!/G}
\newcommand{\xu}{X/\!/U}

\newcommand{\cplusr}{(\CC^+)^r}
\newcommand{\slrplus}{SL(r+1;\CC)}
\newcommand{\glr}{GL(r;\CC)}
\newcommand{\prrplus}{{\PP(\CC \oplus ((\CC^r)^* \otimes \CC^{r+1}))}}
\newcommand{\prr}{{\PP(\CC \oplus ((\CC^r)^* \otimes \CC^{r}))}}
\newcommand{\ximp}{X_{{\rm impl}}}
\newcommand{\tkimp}{(T^*K)_{{\rm impl}}}
\nc{\lieq}{{\mathfrak q}}
\nc{\liez}{{\mathfrak z}}
\nc{\lieqs}{{\lieq}^*}
\nc{\lieg}{{\mathfrak g}}
\nc{\liegs}{{\lieg}^*}
\nc{\liep}{{\mathfrak p}}
\nc{\lieps}{{\liep}^*}
\newcommand{\ximpq}{X_{{\rm impl}}^{{K,K^{(P)}}}}
\newcommand{\tkimpq}{(T^*K)_{{\rm impl}}^{{K,K^{(P)}}}}
\newcommand{\tcone}{\liets_{(P)+}}

\def\a{\alpha}

\def\e{\epsilon}
\def\z{\zeta}

\def\l{\lambda}

\def\x{\xi}

\def\s{\sigma}

\title{Symplectic implosion and non-reductive quotients}

\author{Frances Kirwan \\Mathematical Institute, Oxford OX1 3BJ, UK}

\begin{document}

\maketitle

\section{Introduction}

There is a close relationship between Mumford's geometric invariant theory (GIT) in (complex)
algebraic geometry and the process of reduction in symplectic geometry.
GIT was developed to construct
quotients of algebraic varieties by 
reductive group actions
 and thus to construct and study moduli 
spaces 
 \cite{GIT,New}.  When a moduli space (or a compactification
of a moduli space) over $\CC$ can be constructed as a GIT quotient of a
complex projective variety by the action of a complex reductive group $G$, then 
it can be identified with a symplectic reduction by a maximal compact 
subgroup $K$ of $G$ and techniques
from symplectic geometry can be used to study its topology (for example \cite{AB,JK,JKKW,K,K2,K3,K4}).
 Many moduli spaces arise as quotients of algebraic
group actions, but the groups concerned are not necessarily reductive, so that
classical GIT does not apply and different methods need to be used
to construct the quotients (cf. e.g. \cite{MK,Kollar}). Nonetheless, in
suitable situations GIT can be generalised to allow us to construct
GIT-like quotients (and compactified quotients) for these actions
\cite{DK,DK2,KPEN}. This paper describes some ways in which such non-reductive
compactified quotients can be studied using symplectic techniques closely
related to the \lq symplectic implosion' construction of Guillemin, Jeffrey and
Sjamaar \cite{GJS}.

More precisely, suppose that $U$ is a maximal unipotent subgroup of a complex reductive group
$G$ acting linearly (with respect to an ample line bundle $L$) on a complex projective variety 
$X$, and suppose that the linear action of $U$ on $X$ 
extends to a linear action of $G$. Then
the ring of invariants $\bigoplus_{k \geq 0} H^{0}(X,L^{\otimes k})^U$ is finitely generated and the
 enveloping quotient $\xu$ (in the sense of \cite{DK}) is the projective variety
$\Proj(\bigoplus_{k \geq 0} H^{0}(X,L^{\otimes k})^U)$ associated to the ring of invariants.
Moreover if $K$ is a maximal compact subgroup of $G$, and $X$ is given a suitable $K$-invariant K\"{a}hler
form, then $\xu$
can be identified with the  imploded cross-section $\ximp$ of $X$ by $K$ in the sense of the symplectic implosion construction of Guillemin, Jeffrey
and Sjamaar \cite{GJS}. Note that here $U$ is the unipotent radical of a Borel subgroup of $G$.
The aim of this paper is to 
generalise symplectic implosion to give a symplectic construction for GIT-like 
(compactified) quotients by the unipotent radical $U$
of any parabolic subgroup $P$ of a complex reductive group $G$, 
when the action extends to an action of $G$. Hence we obtain
 a \lq moment
map' description of such compactified 
quotients 
of projective varieties by unipotent radicals of parabolics
which is analogous to the description of a reductive GIT quotient
$Y/\!/G$ as a symplectic quotient $\mu^{-1}(0)/K$ where 
$K$ is a maximal compact subgroup of $G$ 
 and $\mu$ 
 is a moment map.

The layout of the paper is as follows. $\S$2 reviews classical GIT and
its relationship with symplectic geometry, while $\S$3 reviews symplectic implosion from \cite{GJS} and 
extends its construction to cover quotients by unipotent radicals of parabolics.
$\S$4 gives a brief description of the results of
\cite{DK,KPEN} on non-reductive actions  and the construction of 
compactified quotients (more details and a much more
leisurely introduction to non-reductive GIT can be found in \cite{DK})
and finally relates them to symplectic implosion. A simple example when 
$G=SL(2;\CC)$ is worked out in detail at the very end of the paper in
Example \ref{lastex}.

\newpage

\subsection{Index of notation}

Notation is introduced in this paper as follows:

\bigskip

$\mu, \liek, K_\z$ \hfill $\S$2.1

$\hat{\calo}_L(X),\xg$ \hfill $\S$2.2

$\ximp, T, \liet, \liets_+, W, [K_\z,K_\z], \Sigma, B, \umax, \overline{G/\umax}^{{\rm aff}}, \Lambda, \Lambda^*_+, V_\lambda, \Pi, \iota, w_0 $ \hfill $\S$3.1

$B^{{\rm op}}, U^{{\rm op}}_{{\rm max}}, V_\lambda^{(T)}, v_\varpi, \calf, \alpha^\vee, S, \widetilde{\ximp} $ \hfill $\S$3.1

$U, P, L^{(P)}, K^{(P)}, S_P, R^+, R(S_P), Q^{(P)}, \liek^{(P)}, \liez^{(P)},
 \overline{G/U}^{{\rm aff}}\!,\, \xu, E^{(P)}  $ \hfill $\S$3.2

$ V^{(P)}_\varpi, V^{K^{(P)}}_\varpi, v^{(P)}_\varpi, v^{(P)}_{\varpi,\lambda}, \pi^{K^{(P)}}, \tcone, \calf^{(P)}, \liek^{(P)*}_+,
K_\z(P), v_\sigma^{(P)}, \ximpq  $ \hfill $\S$3.2

$\widetilde{\ximpq}, \widetilde{G/U}^{{\rm aff}} $ \hfill $\S$3.3

$X^{ss}, X^s, X^{nss}, X^{ns}, \xu, \overline{G \times_U X}, X^{\bar{s}},
X^{\bar{ss}}  $ \hfill $\S$4.1

$\hat{U}, \widehat{\xu}, \hat{L}_\epsilon = \hat{L}_\epsilon^{(N)}, \calx,\tilde{\calx}, \widetilde{\xu}  $ \hfill $\S$4.2 

\section{Symplectic reduction and geometric invariant theory}

The GIT quotient construction in complex algebraic
geometry is closely related to the process of 
reduction in symplectic geometry. 

\subsection{Symplectic reduction}

Suppose that a compact, connected 
Lie group $K$ with Lie algebra ${\liek}$ acts smoothly
on a symplectic manifold
$X$ and preserves the symplectic form $\omega$. 
 Let us denote
the vector field on $X$ defined by
the infinitesimal action of $a\in {\liek}$ by
$x\mapsto a_x.$
Recall that a moment map for the action of $K$ on $X$ is then a smooth map
$\mu :X\rightarrow {\liek}^{\ast}$
which satisfies
$$d\mu(x)(\xi).a=\omega_x(\xi,a_x)$$
for all $x\in X$, $\xi\in T_xX$ and $a\in {\liek}$. Equivalently, 
if $\mu_a:X \to {\RR}$ denotes the component 
of $\mu$ along
$a\in {\liek}$ defined for all $x\in X$ by the pairing
$\mu_a(x)=\mu(x).a$
between $\mu(x) \in {\liek}^{\ast}$ and
$a \in {\liek}$, then $\mu_a$ is a Hamiltonian function                       
for the vector field on $X$ induced by
$a$. We shall assume that any moment map
$\mu :X\rightarrow {\liek}^{\ast}$ is $K$-equivariant with respect to the
given action of $K$ on $X$ and the coadjoint action of $K$ on $\lieks$.
If the stabiliser $K_{\zeta}$ of $\zeta\in {\liek}^{\ast}$
acts freely on $\mu^{-1}(\zeta)$ then $\mu^{-1}(\zeta)$ is
a submanifold of $X$ and the symplectic form $\omega$ induces a
symplectic structure on the quotient $\mu^{-1}(\zeta)/K_{\zeta}$ which is the Marsden-Weinstein
reduction, or symplectic reduction, at $\zeta$ of the action
of $K$ on $X$. The quotient
$\mu^{-1}(\zeta)/K_{\zeta}$ also inherits a symplectic structure
when the action of $K_{\zeta}$
on $\mu^{-1}(\zeta)$ is not free, but in this case it is likely
to have singularities (although these will only
be orbifold singularities if $\zeta$ is a regular value
of $\mu$, or equivalently if $K_\zeta$ acts on $\mu^{-1}(\zeta)$
with finite stabilisers). The case when $\zeta = 0$ is of particular
importance; $\mu^{-1}(0)/K$ is often called the symplectic
quotient of $X$ by the action of $K$.

Now let 
$X$ be a nonsingular connected complex projective variety
embedded in complex projective space $\PP^n$, and let $G$
be a complex Lie group acting on $X$ via a complex linear
representation $\rho:G\rightarrow GL(n+1;\CC)$. 
By an appropriate choice
of coordinates on $\PP^n$ we may assume that $\rho$ maps 
a maximal compact subgroup $K$
of $G$ into the unitary group $U(n+1)$. Then the Fubini-Study form $\omega$ on $\PP^n$ restricts to
a $K$-invariant K\"{a}hler form on $X$, and there is a moment map
$\mu :X\rightarrow {\liek}^{\ast}$ defined (up to multiplication by a 
constant scalar factor depending on the convention chosen for
the normalisation of the Fubini-Study form) by
\begin{equation} \mu(x).a = \frac{\overline{\hat{x}}^{t}\rho_{\ast}(a)\hat{x}}
{2\pi i|\!|\hat{x}|\!|^2} \label{mmap} \end{equation}
for all $a\in {\liek}$, where $\hat{x}\in {\CC}^{n+1}-\{0\}$ is a representative
vector for $x\in \PP^n$ and the representation $\rho:K \to U(n+1)$ induces
$\rho_*: \liek \to {\lieu}(n+1)$ and dually $\rho^*:{\lieu}(n+1)^* \to \lieks$.

In this situation there are two possible quotient
constructions: the symplectic reduction 
$\mu^{-1}(0)/K$  in symplectic geometry and 
the GIT quotient $X/\!/G$  in algebraic geometry described below. In fact these give us
the same space, at least up to homeomorphism (and
diffeomorphism away from the singularities).

\subsection{Mumford's geometric invariant theory}

Let $X$ be a complex projective variety and let $G$ be a complex
reductive group acting on $X$. 
Recall that over $\CC$ a linear algebraic group $G$ is reductive if and only if it is
the complexification of a maximal compact subgroup $K$. The simplest non-trivial
example is the
complexification $\CC^*$ of the circle $S^1$, and more generally
$GL(n;\CC)$ is the complexification of the unitary group
$U(n)$ and thus is reductive. In contrast the additive group of complex
numbers $\CC^+$ has no nontrivial compact subgroups and so is not reductive;
the same is true of any complex linear algebraic group $U$ which is unipotent (that is,
$U$ is isomorphic to a closed subgroup of the group of strictly upper triangular matrices in $GL(n;\CC)$ for some $n$). In some sense reductive and unipotent groups sit at the opposite extremes of a spectrum, and any linear algebraic group
$H$ has a unique maximal unipotent normal subgroup $U$ (its unipotent radical)
such that the quotient group $H/U$ is reductive.

Geometric invariant theory needs an extra ingredient in addition to the
action of $G$ on $X$, which is
 a  {linearisation} of the action; that is, a
line bundle $L$ on $X$ and a lift of the action of $G$ to $L$. The line bundle
$L$ is usually taken to be ample, and then
very little generality is lost by assuming that for some
projective embedding 
$$X \subseteq \PP^n$$ 
the action
of $G$ on $X$ extends to an action on $\PP^n$ given by a
representation 
$$\rho:G\rightarrow GL(n+1),$$
and taking for $L$ the hyperplane line bundle $\calo_{\PP^n}(1)$ on $\PP^n$.

A {\em categorical quotient} of a variety $X$ under an
action of  $G$ is a $G$-invariant morphism $\phi:X \to Y$ from $X$ to a variety $Y$ 
such that any other $G$-invariant morphism $\tilde{\phi}: X \to
\tilde{Y}$ factors as $\tilde{\phi} = \chi \circ \phi$ for a unique
morphism $\chi:Y \to \tilde{Y}$
\cite[Chapter 2, $\S$4]{New}. An {\em orbit space} for the action
is a categorical quotient $\phi:X \to Y$ such that each fibre $\phi^{-1}(y)$ is
a single $G$-orbit, and a {\em geometric
quotient} is an orbit space $\phi:X \to Y$ which is an affine morphism such that

(i) if $U$ is open in $Y$ then
$$ \phi^*: {\calo}(U) \to {\calo}(\phi^{-1}(U)) $$
induces an isomorphism of ${\calo}(U)$ onto ${\calo}(\phi^{-1}(U))^G$, and

(ii) if $W_1$ and $W_2$ are disjoint closed $G$-invariant
subvarieties of $X$ then their images $\phi(W_1)$ and $\phi(W_2)$ in
$Y$ are disjoint closed subvarieties of $Y$.

When $G$ acts linearly on $X$ as above there is an induced action of $G$ on the
homogeneous coordinate ring 
\begin{equation} {\hat{\calo}}_L(X) = \bigoplus_{k \geq 0} H^{0}(X, L^{\otimes k}) \cong
\CC[x_0,...,x_n]/\cali_X \end{equation}
where $\cali_X$ is the ideal in $\CC[x_0,...,x_n]$
generated by the homogeneous polynomials vanishing on $X$.
The
subring ${\hat{\calo}}_L(X)^G$ consisting of the elements of ${\hat{\calo}}_L(X)$
left invariant by $G$ is a finitely generated graded complex algebra
because $G$ is reductive, and
so we can define the GIT quotient $X/\!/G$ to be the projective variety  
$\Proj({\hat{\calo}}_L(X)^G)$ associated to ${\hat{\calo}}_L(X)^G$ \cite{GIT}. The inclusion of
${\hat{\calo}}_L(X)^G$ in ${\hat{\calo}}_L(X)$ determines a { rational} map $q$ from $X$
to $X/\!/G$, but in general there will be points of $X \subseteq \PP^n$ where
every $G$-invariant polynomial vanishes and so this map will not be
well-defined everywhere on $X$. Hence
we define the set $X^{ss}$ of {\em semistable} points in $X$ to
be the set of those $x \in X$ for which there exists
some $f \in {\hat{\calo}}_L(X)^G$ not vanishing at $x$, and then the rational
map $q$ restricts to a surjective $G$-invariant morphism from
the open subset $X^{ss}$ of $X$ to the quotient variety $X/\!/G$,
which is a categorical quotient for the action of $G$ on $X^{ss}$. This restriction
$q:X^{ss} \to X/\!/G$ is not necessarily an orbit
space: when $x$ and $y$ are semistable points of $X$ we have
$q(x)=q(y)$ if and only if the closures $\overline{O_G(x)}$ and
$\overline{O_G(y)}$ of the $G$-orbits of $x$ and $y$ meet in
$X^{ss}$. Topologically $X/\!/G$ is the quotient of $X^{ss}$ by the
equivalence relation $\sim$ such that if $x$ and $y$ lie in $X^{ss}$ 
then $x \sim y$ if and only if $\overline{O_G(x)}$ and
$\overline{O_G(y)}$ meet in
$X^{ss}$.

A {\em stable} point of $X$ (\lq properly stable' in
the terminology of \cite{GIT})
is a point $x$
of $X^{ss}$ with a $G$-invariant neighbourhood in $X^{ss}$ such that every $G$-orbit 
in this neighbourhood
is closed in $X^{ss}$ and has dimension
 $\dim G$. 
If $U$ is any $G$-invariant open
subset of the set $X^s$ of stable points of $X$, then $q(U)$ is
an open subset of $X/\!/G$ and the restriction
$q|_U :U \to q(U)$ of $q$ to $U$ is an orbit space
for the action of $G$ on $U$,
so that it makes sense to write $U/G$ for $q(U)$; in fact $U/G$ is
a geometric quotient for the action of $G$ on $U$.
In particular there is a geometric quotient $X^s/G$ for the action
of $G$ on $X^s$, and $X/\!/G$ can
be thought of as a compactification of $X^s/G$.

\begin{equation} \label{gitpicture}
\begin{array}{ccccc}
   X^s & \subseteq & X^{ss} & \subseteq & X \\
       & {\rm open} &  & {\rm open} & \\
  \Big\downarrow & & \Big\downarrow & & \\
  &  &  &  &  \\
  X^s/G & \subseteq &   X/\!/G = X^{ss} / \sim & & \\
       & {\rm open} & & & 
\end{array} \end{equation}

\bigskip

\begin{rem} $X^s, X^{ss}$ and $X/\!/G$ are unaltered if  for any $k > 0$ the
line bundle $L$ is replaced by $L^{\otimes k}$
with the induced action of $G$, so it is sometimes convenient to allow
{\em fractional} linearisations $L^{\otimes {\ell}/{m}}$.
\end{rem}

The subsets $X^{ss}$ and $X^s$ of $X$ are characterised by the following properties (see
Chapter 2 of \cite{GIT} or \cite{New}). 

\begin{prop} (Hilbert-Mumford criteria)
\label{sss} (i) A point $x \in X$ is semistable (respectively
stable) for the action of $G$ on $X$ if and only if for every
$g\in G$ the point $gx$ is semistable (respectively
stable) for the action of a fixed maximal (complex) torus of $G$.

\noindent (ii) A point $x \in X$ with homogeneous coordinates $[x_0:\ldots:x_n]$
in some coordinate system on $\PP^n$
is semistable (respectively stable) for the action of a maximal (complex)
torus of $G$ acting diagonally on $\PP^n$ with
weights $\a_0, \ldots, \a_n$ if and only if the convex hull
$$\conv \{\a_i :x_i \neq 0\}$$
contains $0$ (respectively contains $0$ in its interior).
\end{prop}

The GIT quotient $\xg$ is homeomorphic to the symplectic quotient
$\mu^{-1}(0)/K$, and the subsets $X^{ss}$ and $X^s$ of $X$ can be described using the moment map
$\mu$ at (\ref{mmap}) above.
More precisely \cite{K},
any $x\in X$ is semistable if and only if
the closure of its $G$-orbit meets $\mu^{-1}(0)$, 
while $x$ is stable if and only if its $G$-orbit
meets
$$\mu^{-1}(0)_{\rm reg} = \{ x \in \mu^{-1}(0)\  | \ d\mu(x):T_xX \to \lieks
\mbox{ is surjective} \},$$
and
the inclusions of $\mu^{-1}(0)$ into $X^{ss}$ and
of $\mu^{-1}(0)_{\rm reg}$ into $X^{s}$
induce 
homeomorphisms
$$\mu^{-1}(0)/K\rightarrow X/\!/G$$
and
$$\mu^{-1}(0)_{\rm reg} \to X^{s}/G.$$
Thus the moment map picks out a unique
$K$-orbit in each stable $G$-orbit, and also in each
equivalence class of strictly semistable
$G$-orbits, where $x$ and $y$ in $X^{ss}$
are equivalent if the closures of their
$G$-orbits meet in $X^{ss}$ (that is,
if their images under the natural
surjection $q:X^{ss} \to X/\!/G$ agree).

\begin{rem}
It follows from the formula (\ref{mmap}) that if we change
the linearisation of the $G$-action of $X$ by multiplying
by a character $\chi:G \to \CC^*$ of $G$, then the
moment map is modified by the addition of a central
constant $c_\chi$ in $\lieks$, which we can identify with 
the restriction to $\liek$ of the derivative of $\chi$.
\end{rem}

\begin{ex} Let $G=SL(2;\CC)$ act on $X=(\PP^1)^4$ via M\"{o}bius
transformations and let $K$ be the maximal compact subgroup $SU(2)$ of $G$.
If we identify $\PP^1$ with the unit sphere $S^2$ in $\RR^3$
then there is a moment map $$\mu:X = (S^2)^4 \to \lieks \cong \RR^3$$ given by
$\mu(x_1,x_2,x_3,x_4) = x_1 + x_2 + x_3 + x_4$. Thus
$\mu^{-1}(0)$ consists of configurations of 4 points on $S^2$ which are balanced
in the sense that their centre of gravity lies at the origin,
while $\mu^{-1}(0) \setminus \mu^{-1}(0)_{\rm reg}$ consists of the
configurations in which two points coincide at some $p \in S^2$ and the
other two points coincide at the antipodal point $-p$.
The open subset
$$X^s = \{ (x_1,x_2,x_3,x_4)\in (\PP^1)^4: x_1,x_2,x_3,x_4 \mbox{ distinct} \}$$
of $X=(\PP^1)^4$ has a geometric quotient which, using the cross-ratio, can be identified 
with
$$\PP^1 - \{0,1,\infty\}$$
and this in turn can be identified with $\mu^{-1}(0)_{\rm reg}/K$. In addition
$$X^{ss} = \{ (x_1,x_2,x_3,x_4)\in (\PP^1)^4: \mbox{ at most two of } x_1,x_2,x_3,x_4 \mbox{ coincide} \}$$
has a categorical quotient $\xg \cong X^{ss}/\! \sim \  \cong \PP^1$ in which the points $0,1,\infty$ each represent three strictly 
semistable $G$-orbits in $X$: one $G$-orbit consisting of 
configurations in which two points $x_i$ and $x_j$ coincide at some $p \in \PP^1$ and the
other two points $x_k$ and $x_m$ coincide at a distinct point $q \in \PP^1$, a second
consisting of 
configurations in which $x_i$ and $x_j$ coincide at some $p \in \PP^1$ and the
other two points $x_k$ and $x_m$ are distinct from each other and from $p$, and
the third  consisting of 
configurations in which $x_k$ and $x_m$ coincide at a some point $q \in \PP^1$
while $x_i$ and $x_j$ are distinct from each other and from $q$. The first of these
orbits is closed in $X^{ss}$ and lies in the closure of each of the other two orbits.

\end{ex}

\section{Symplectic implosion and quotients by non-reductive groups}


 Ways in which classical GIT might be generalised to actions of non-reductive affine
 algebraic  groups on algebraic varieties were studied in \cite{DK} (see also \cite{KPEN})
 building on earlier work such as \cite{F2,F1,GP1,GP2,W}.
 Every affine algebraic group $H$ has a unipotent radical $U \unlhd
 H$ such that $H/U$ is reductive, so we can concentrate on unipotent actions. It is
 shown in \cite{DK} that when a unipotent group $U$ acts linearly (with respect to an ample
 line bundle $L$) on a complex projective variety $X$, then $X$ has invariant
 open subsets $X^s \subseteq X^{ss}$, consisting of the \lq stable' and \lq semistable'
 points for the action, such that $X^s$ has a 
 geometric quotient $X^s/U$ and $X^{ss}$ has a 
 canonical \lq enveloping quotient'  $X^{ss} \to X/\!/U$ which restricts to $X^s \to
 X^s/U$ where $X^s/U$ is an open subset of $X/\!/U$. However, in contrast to the reductive case, the natural map
from $X^{ss}$ to $X/\!/U$ is not necessarily surjective, and indeed
 its image is not necessarily a subvariety of $X/\!/U$, so this does not in general
 give us a categorical quotient of $X^{ss}$. Furthermore $X/\!/U$ is in general only 
 quasi-projective, not projective, though when the ring of invariants
 $\hat{\calo}_L(X)^U = \bigoplus_{k \geq 0} H^{0}(X, L^{\otimes k})^U$
is finitely generated as a $\CC$-algebra then
$X/\!/U$ is the projective variety $\Proj({\hat{\calo}}_L(X)^U)$.

In order to obtain a compactification $\overline{\xu}$ 
of the enveloping quotient
$X/\!/U$ when the ring of invariants ${\hat{\calo}}_L(X)^U$ is not
finitely generated, and to understand its geometry even
when $\xu = \overline{\xu}$ is itself projective, we can
 transfer the problem of constructing a quotient for the
$U$-action to the construction of a quotient for an action of a reductive
group $G$ which contains $U$ as a subgroup, by finding a reductive
envelope. This is a projective completion 
$$\overline{G \times_U X}$$
of the quasi-projective variety $G \times_U X$ (which is the quotient of $G \times X$ by the free action of $U$ acting diagonally on the left on $X$
and by right multiplication on $G$), with a linear $G$-action on 
$\overline{G \times_U X}$ extending the induced $G$-action
on $G \times_U X$, such that the $U$-invariants on
$X$ lying in a suitable  set (see Definition \ref{defn:separ} below)
 extend to $G$-invariants on $\overline{G \times_U X}$.
If  the linearisation on
$\overline{G \times_U X}$ is ample, then the classical GIT quotient
$$\overline{G \times_U X}/\!/G$$
is a compactification
$\overline{X/\!/U}$ of $X/\!/U$, and hence also of its
open subset $X^s/U$ if $X^s \neq \emptyset$. Moreover 
if $X^{\bar{s}}$ and $X^{\bar{ss}}$ denote the
open subsets of $X$ consisting of points of $X$ which are stable
and semistable for the $G$-action on
$\overline{G \times_U X}$ under the inclusion 
$$X \hookrightarrow G \times_U X \hookrightarrow \overline{G \times_U X}$$
then
$$X^{\bar{s}} \subseteq X^s \subseteq X^{ss} \subseteq X^{\bar{ss}}.$$
Note however that $X^{\bar{s}}, X^{\bar{ss}}$ and $\overline{X/\!/U}$
depend in general on the choice of reductive envelope
$\overline{G \times_U X}$ with its linear $G$-action, whereas $X^{{s}}, X^{{ss}}$ and ${X/\!/U}$
depend only on the linear action of $U$ on $X$.

Just as GIT quotients by complex reductive groups are closely
related to symplectic reduction, so quotients by suitable unipotent
groups (in particular maximal unipotent subgroups of complex reductive
groups)
are closely related to the construction 
called symplectic implosion \cite{GJS} which we will discuss below. 

\subsection{Symplectic implosion for a maximal unipotent subgroup}

Let $(X,\omega)$ be a symplectic manifold on which a compact connected Lie group $K$
acts with a moment map $\mu:X \to \lieks$ where $\liek$ is the Lie algebra of $K$. Let
us choose an invariant inner product on $\liek$ and use it to identify $\lieks$
with $\liek$. Let $T$ be a maximal torus of $K$ with Lie algebra $\liet \subseteq \liek$
and Weyl group $W = N_K(T)/T$, and let $\liets_+ \cong \liets/W \cong \lieks/{\rm Ad}^*(K)$
be a positive Weyl chamber in $\lieks$. The {\em imploded cross-section} \cite{GJS}
of $X$ is then
\begin{equation} X_{{\rm impl}} = \mu^{-1}(\liets_+)/\approx \end{equation}
where $x \approx y$ if and only if $\mu(x) = \mu(y) = \zeta \in \liets_+$ and
$x = ky$ for some $k \in [K_\zeta,K_\zeta]$. Here $K_\zeta$ denotes the stabiliser
$K_\zeta = \{ k \in K: ({\rm Ad}^*k) \zeta = \zeta \}$ of $\zeta$ under the co-adjoint action of
$K$ on $\lieks$, and $[K_\zeta,K_\zeta]$ is its commutator subgroup. If $\Sigma$ is the 
set of faces of $\liets_+$ then 
\begin{equation}  X_{{\rm impl}}  \ = \  \coprod_{\sigma \in \Sigma} \frac{\mu^{-1}(\sigma)}{[K_\sigma,K_\sigma]}
\ =  \ \mu^{-1}((\liets_+)^\circ) \ \ \ \sqcup
 \coprod_{\begin{array}{c}\sigma \in \Sigma\\  \sigma \neq (\liets_+)^\circ \end{array}
 } \frac{\mu^{-1}(\sigma)}{[K_\sigma,K_\sigma]} \end{equation}
where $K_\sigma = K_\zeta$ for any $\zeta \in \sigma$. The topology on $\ximp$
is the quotient topology induced from $\mu^{-1}(\liets_+)$, and $\ximp$ also 
inherits a symplectic structure. More precisely, it is stratified by the locally
closed subsets $\mu^{-1}(\sigma)/[K_\sigma,K_\sigma]$, each of which is the symplectic
reduction by the action of $[K_\sigma,K_\sigma]$  of a locally closed symplectic
submanifold
$$X_\sigma = K_\sigma \mu^{-1}( \bigcup_{\tau \in \Sigma, \bar{\tau} \supseteq \sigma} \tau)$$
of $X$ (and locally near every point $\ximp$ can be identified symplectically with the product of the
stratum and a suitable cone in the normal direction). The induced action of $T$ on 
$\ximp$ preserves this symplectic structure and has a moment map
$$\mu_{\ximp}:\ximp \to \liets_+ \subseteq \liets$$
inherited from the restriction of $\mu$ to $\mu^{-1}(\liets_+)$. 
 If $\zeta \in \liets_+$ the symplectic reduction of $\ximp$ at
 $\zeta$ for this action of $T$ is the symplectic reduction of $X$ at $\zeta$ for the action
 of $K$:
\begin{equation} \frac{\mu_{\ximp}^{-1}(\zeta)}{T} = \frac{\mu^{-1}(\zeta)}{T.[K_\zeta,K_\zeta]} = \frac{\mu^{-1}(\zeta)}{K_\zeta}.
\end{equation}
The {\em universal imploded cross-section} is the imploded cross-section
\begin{equation} (T^*K)_{{\rm impl}} = K \times \liets_+ / \approx \end{equation}
of the cotangent bundle $T^*K \cong K \times \lieks$ with respect to the $K$-action induced from the 
right action of $K$ on itself; it inherits an action of $K \times T$ from the left action
of $K$ on itself and the right action of $T$ on $K$. Any other imploded cross-section
$\ximp$ can be constructed as the symplectic quotient of the product $X \times \tkimp$
by the diagonal action of $K$ (\cite{GJS} Theorem 4.9).

In fact $\tkimp$ is always a complex affine variety and its symplectic structure is
given by a K\"{a}hler form. Indeed, let $G=K_c$ be the complexification of $K$ and
let $B$ be a Borel subgroup of $G$ with $G = KB$ and $K \cap B = T$. If $\umax \leq B$
is the unipotent radical of $B$ (and hence a maximal unipotent subgroup of $G$), then
$\umax$ is a Grosshans subgroup of $G$ \cite{Grosshans}: that is, the quasi-affine variety $G/\umax$
can be embedded as an open subset of an affine variety in such a way that its
complement has (complex) codimension at least two. This means that the ring 
of invariants $\calo(G)^\umax$ is 
finitely generated (see for example \cite{Grosshans}), and by
\cite{GJS} Proposition 6.8 there is a natural $K \times T$-equivariant identification 
$$\tkimp \cong {\rm Spec}(\calo(G)^\umax)$$
of the canonical affine completion ${\rm Spec}(\calo(G)^\umax)$
of $G/\umax$ with $\tkimp$. It follows that if $X$ is a complex projective
variety on which $G$ acts linearly with respect to a very ample
line bundle $L$, and $\omega$ is an associated $K$-invariant K\"{a}hler form on
$X$, then the symplectic quotient $\ximp$ of $X \times \tkimp$ by $K$
can be identified with the GIT quotient $(X \times {\rm Spec} (\calo(G)^\umax )))/\!/G$.
Moreover
$$\hat{\calo}_L(X)^\umax \cong (\hat{\calo}_L(X) \otimes \calo(G)^\umax)^G$$
is finitely generated, and if we define the GIT quotient 
$X/\!/\umax$ to be the projective variety ${\rm Proj}(\hat{\calo}_L(X)^\umax)$
associated to the ring of invariants $\hat{\calo}_L(X)^\umax$ then
\begin{equation} X/\!/\umax = {\rm Proj}(\hat{\calo}_L(X)^\umax) \cong (X \times {\rm Spec} (\calo(G)^\umax )))/\!/G
\cong \ximp. \end{equation}

The proof in \cite{GJS} $\S$6 that $\tkimp$ is homeomorphic to the canonical affine completion 
$$ \overline{G/\umax}^{{\rm aff}} = {\rm Spec}(\calo(G)^\umax)$$
of $G/\umax$ runs as follows.
First it is possible to reduce to the case when $K$ is semisimple and simply connected, by regarding
$K$ as the quotient by a finite central subgroup of
$Z(K) \times \widetilde{[K,K]}$ where $Z(K)$ is the centre of $K$
and $\widetilde{[K,K]}$ is the universal cover of the
commutator subgroup $[K,K]$ of $K$.

Following \cite{GJS} $\S$6, if $K$ is a semisimple, connected and simply connected compact group
let $\Lambda = \ker(\exp |_{\liet})$ be the exponential lattice in $\liet$, and let
$\Lambda^* = {\rm Hom}_{\ZZ}(\Lambda,\ZZ)$ be the weight lattice in $\liets$, so that 
$\Lambda^*_+ = \Lambda^* \cap \liets_+$ is the monoid of dominant weights. For $\lambda \in \Lambda^*_+$
let $V_{\lambda}$ be the irreducible $G$-module with highest weight $\lambda$, and let
$$\Pi = \{\varpi_1, \ldots, \varpi_r \}$$
be the set of fundamental weights, which forms a $\ZZ$-basis of $\Lambda^*$ and a minimal set of generators
for $\Lambda^*_+$. 
Recall that $V_\lambda^* = V_{\iota \lambda}$ is the irreducible $G$-module with highest
weight $\iota \lambda$, where $\iota:\liets \to \liets$ is the involution given by
$\iota \l = - w_0 \l$ and $w_0$ denotes the element of the Weyl group $W$ of $G$ such that
$w_0 \umax w_0^{-1} = U_{{\rm max}}^{{ \rm op}}$ is the unipotent radical of the Borel subgroup
$B^{{\rm op}}$ of $G$ which is opposite to $B \geq U$ in the sense that
$B \cap B^{{\rm op}}$ is the complexification $T_c$ of $T$ and $\umax \cap U_{{\rm max}}^{{\rm op}}$
is the identity subgroup.
We have an isomorphism of $G$-modules
\begin{equation} \label{thisisobb}
\calo(G)^\umax \cong \bigoplus_{\lambda \in \Lambda^*_+} V_{\lambda}^* 
\cong \bigoplus_{\lambda \in \Lambda^*_+} V_{\iota \lambda}\end{equation}
where $G$ acts on itself on the left and $\umax$ acts on $G$ on the right. Note that
$T_c$ normalises $\umax$ and this isomorphism (\ref{thisisobb})
becomes an isomorphism of $G \times T_c$-modules if we let $T_c$ act on 
$V_{\lambda}$ with weight $-\lambda$ so that it acts on 
$V_{\lambda}^*$ with weight $\lambda$
(see \cite{Grosshans} $\S$12). 
Equivalently we have an isomorphism of $G \times T_c$-modules
\begin{equation} \label{thisiso2bb}
\calo(G)^\umax \cong \bigoplus_{\lambda \in \Lambda^*_+}  V_\l^{(T)} \otimes V_{\lambda}^* \end{equation}
where $V^{(T)}_\l$ is the irreducible $T_c$-module with  weight $\l$, and
by \cite{Grosshans} Theorem 12.9 this isomorphism extends to an isomorphism of
$G \times G$-modules
\begin{equation} \label{thisiso3bb}
\calo(G) \cong \bigoplus_{\lambda \in \Lambda^*_+} V_{\lambda} \otimes V_\l^*. \end{equation}
In particular the algebra $\calo(G)^\umax$ is generated
 by its finite-dimensional vector subspace
$$\bigoplus_{\varpi \in \Pi} V_{\varpi}^*
\cong \bigoplus_{\varpi \in \Pi} V_{\varpi}^{(T)} \otimes V_{\varpi}^*.$$
The inclusion of this finite-dimensional subspace into $\calo(G)^\umax$ induces a closed
$G \times T_c$-equivariant embedding of $\overline{G/U}_{{\rm max}}^{{\rm aff}} = {\rm Spec}(\calo(G)^\umax)$
into the affine space
$$E = \bigoplus_{\varpi \in \Pi} V_{\varpi} \cong \bigoplus_{\varpi \in \Pi} (V_{ \varpi}^{(T)})^* \otimes
 V_{\varpi}
,$$
sending the identity coset $\umax$ in $ G/\umax \subseteq \overline{G/U}_{{\rm max}}^{{\rm aff}}$
to a sum
$$ \sum_{\varpi \in \Pi} v_{\varpi}$$
of highest weight vectors $v_{\varpi} \in V_{\varpi} \cong (V_{ \varpi}^{(T)})^* \otimes
 V_{\varpi}$. Under this embedding $G/\umax$ is identified with
 $G \, E^{\umax}$ where $E^{\umax}$ is the subspace of $E$ consisting of vectors
 fixed by $\umax$.
 We give $E$ a flat K\"{a}hler
structure $\omega_E$ via the unique $K \times T$-invariant Hermitian inner product on $E$ which satisfies
$|\!| v_{\varpi}|\!| = 1$ for each $\varpi \in \Pi$. Then by \cite{GJS} Proposition 6.8
there is a $K \times T$-equivariant map $\calf:K \times \liets_+ \to E$ defined
on $\liets_+$ by
\begin{equation} \label{calff}
\calf(1,\lambda) = \frac{1}{\sqrt{\pi}} \sum_{j=1}^r \sqrt{\lambda(\alpha_j^{\vee})} v_{\varpi_j},
\end{equation}
where $\alpha^{\vee} = 2\alpha/(\alpha\cdot\alpha)$ and
$$ S = \{ \alpha_1, \ldots ,\alpha_r \}$$ is the set of simple roots corresponding to the
fundamental weights $\{ \varpi_1, \ldots, \varpi_r \}$ (so that
$\varpi_i . {\alpha}_j^\vee = \delta_{ij}$ for $i,j \in \{1,\ldots,r\}$);
moreover $\calf$ induces a homeomorphism from $\tkimp$ to
$\overline{G/\umax}^{{\rm aff}}$ whose restriction to each stratum $\mu^{-1}(\sigma)/[K_\sigma,K_\sigma]$
of $\tkimp$ is a symplectic isomorphism onto its image.

\begin{rem} \label{calffo}
Let $M$ be any compact K\"{a}hler manifold on which the complexified torus $T_c$ acts in such a
way that $T$ preserves the K\"{a}hler structure and has a moment map $\mu_T:M \to \liets$.
In \cite{Atiyah} Theorem 2 Atiyah shows 

(a) that the image $\mu_T(\bar{Y})$ under the torus moment
map $\mu_T$ of the closure $\bar{Y}$ in $M$ of the $T_c$-orbit $Y = T_c m$ of any $m \in M$ is a convex polytope
$\calp$ whose vertices are the images under $\mu_T$ of the connected components of $\bar{Y} \cap M^T$
where $M^T$ is the $T$-fixed point set in $M$, 

(b) that the inverse image in $\bar{Y}$ of each open face of $\calp$ consists of a single
$T_c$-orbit, and

(c) that $\mu_T$ induces a homeomorphism of $\bar{Y}/T$ onto $\calp$.

\noindent In fact Atiyah's proof shows that if $\caly = \exp(i\liet)$ is the orbit of $m \in M$ under the subgroup $\exp(i\liet)$ of $T_c$ then $\mu_T$ restricts to a homeomorphism from $\bar{\caly}$ onto $\calp$, and the
inverse image in $\bar{\caly}$ of each open face of $\calp$ consists of a single $\exp(i\liet)$-orbit.

We can apply this to the compactification $M = \PP(\CC \oplus E)$ of the affine space $E$. The
moment map $\mu_T^E:E \to \liets$ for the $T$-action on $E$ with its chosen flat K\"{a}hler structure
is given (up to multiplication by a positive constant) by
$$\sum_{\varpi} u_{\varpi} \mapsto \sum_{\varpi} |\!| u_{\varpi}|\!|^2 \varpi$$
when $u_{\varpi} \in V_{\varpi}$ for $\varpi \in \Pi $, while the 
moment map $\mu_T^{\PP(\CC\oplus E)}:\PP(\CC\oplus E) \to \liets$ for the $T$-action on $\PP(\CC\oplus E)$ with 
the induced Fubini-Study K\"{a}hler structure
is given (up to multiplication by a positive constant) by
$$[z: \sum_{\varpi \in \Pi} u_{\varpi}] \mapsto \frac{\sum_{\varpi} |\!| u_{\varpi}|\!|^2 \varpi}
{|z|^2 + \sum_{\varpi \in \Pi} |\!| u_{\varpi}|\!|^2}$$
when $z \in \CC$ and $u_{\varpi} \in V_{\varpi}$ for $\varpi \in \Pi $
are not all zero.
Comparing these two moment maps on $E$ (regarded as as open subset of $\PP(\CC\oplus E)$ in the usual way)
we see that the image under $\mu_T^E$ of the closure $\bar{\caly}$ in $E$ of the $\exp(i\liet)$-orbit $\caly$
in $E$ of the vector $\sum_{\varpi \in \Pi} v_{\varpi}$ corresponding to the identity coset $\umax$ in $G/\umax$ is 
the cone in $\liets$ spanned by the half-lines $\RR_+\varpi$ for $\varpi \in \Pi$, which is
of course the positive Weyl chamber $\liets_+$. We also find that the
restriction 
\begin{equation} \label{muTE} \mu_T^E|_{\bar{\caly}}: \bar{\caly} \to \liets_+ \end{equation}
is a homeomorphism, and it is easy to check that the map $\calf:\liets_+ \to E$ of \cite{GJS} Proposition
6.8 defined at (\ref{calff}) above can be identified with the composition of 
the inverse $(\mu_T^E|_{\bar{\caly}})^{-1}:  \liets_+ \to \bar{\caly}$ of (\ref{muTE}) and the
inclusion of $\bar{\caly}$ in $E$. From this it can be deduced that its $K \times T$-equivariant
extension $\calf: K \times \liets_+ \to E$ induces a bijection from $\tkimp$ onto the
closure $\overline{G/\umax}^{{\rm aff}}$ of $G(\sum_{\varpi \in \Pi} v_{\varpi}) \cong G/\umax$ in $E$
using

(i) the Iwasawa decomposition
$$G = K \ \exp(i\liet)\  \umax$$
of $G$ which tells us that $\overline{G/\umax}^{{\rm aff}} = K \bar{\caly} = \calf(K \times \liets_+)$, and

(ii) Lemma 6.2 of \cite{GJS}, which shows that for each face $\sigma$ of 
$\liets_+$ the stabiliser
in $K$ of $$\sum_{\varpi \in \sigma} v_{\varpi}$$ is $[K_\sigma, K_\sigma ]$.
\end{rem}

Guillemin, Jeffrey and Sjamaar also construct a $K\times T$-equivariant desingularisation
$\widetilde{\tkimp}$ for the universal imploded cross-section $\tkimp \cong \overline{G/\umax}^{{\rm aff}}$
and a partial desingularisation $\widetilde{\ximp}$ for $\ximp$. In \cite{GJS} $\S$7 they show that 
if the action of $K$ on $X$ 
has principal face the interior $(\liets_+)^\circ$ of $\liets_+$ (where
the principal face is the minimal open face $\s$ of $\liets_+$ such that
$\mu(X) \cap \liets_+$ is contained in $\bar{\s}$), then $\widetilde{\ximp}$ can be identified with
the symplectic quotient of $X \times \widetilde{\tkimp}$ by the induced action of $K$ (and they
observe without proof that the same is true for any principal face). 
Moreover $\widetilde{\tkimp}$ can be identified as a Hamiltonian $K$-manifold with the
homogeneous complex vector bundle
\begin{equation} \widetilde{G/\umax}^{{\rm aff}} = G \times_B E^{\umax} \end{equation}
 over the flag manifold $G/B$, where the restriction 
to $G \times E^{\umax}$ of the multiplication map $G \times E \to E$ induces a birational $G$-equivariant
morphism 
$$p_{\umax}: \widetilde{G/\umax}^{{\rm aff}} \to \overline{G/\umax}^{{\rm aff}} = \tkimp \subseteq E.$$
Note that the fixed point set $E^{\umax}$ of $\umax$ in $E$ is the closure in $E$ of the $T_c$-orbit of $\sum_{\varpi \in \Pi} v_{\varpi}$.
If $\lambda_0 \in \liets$ is regular dominant and $\e >0$ is sufficiently close to 0, and if $\omega_0$ is the K\"{a}hler form
on $G/B$ given by regarding $G/B$ as the coadjoint $K$-orbit through $\e \lambda_0$, then $p_{\umax}^* \omega_E + q^* \omega_0$ is a K\"{a}hler form on $\widetilde{G/\umax}^{{\rm aff}}$ where $q: G \times_B E \to G/B$ is the
projection.

It is also shown in \cite{GJS} $\S$7 that 
the partial desingularisation $\widetilde{\ximp}$ can alternatively be obtained from
$\ximp$ via a symplectic cut with respect to the $T$-action and the polyhedral
cone $\e \lambda_0 + \overline{\tau}$ where $\tau$ is the principal face of $X$ and $\lambda_0 \in \tau$ and $\e >0$ is sufficiently close to 0;
that is, $\widetilde{\ximp}$ is the symplectic reduction at $\e \lambda_0$ for the diagonal
$T$-action on the product of $\ximp$ and the symplectic toric manifold
associated to the polyhedron $-\overline{\tau}$ (see \cite{Lerman,meinsj}).

\subsection{Symplectic implosion for the unipotent radical of a parabolic subgroup}

Now suppose that $U$ is the unipotent radical of a parabolic subgroup $P$ of the complex reductive
group $G$. Recall (see e.g. \cite{borel, Springer}) that a parabolic subgroup of $G$ is a closed subgroup which contains some Borel subgroup, and its
unipotent radical is its unique maximal normal unipotent subgroup;
thus by replacing $P$ with a suitable conjugate in $G$ if necessary,
we can assume that $P$ contains the Borel subgroup $B$ of $G$ and $U \leq \umax$. Then $P = 
U {L^{(P)}} \cong U \rtimes {L^{(P)}}$, where
the Levi subgroup ${L^{(P)}}$ of $P$ contains the complex maximal torus $T_c$ of $G$,
and we can assume in addition that ${L^{(P)}}$ is the complexification of its intersection
$${K^{(P)}} = {L^{(P)}} \cap K = P \cap K$$
with $K$. For some subset $S_P$ of the set $S$ of simple roots, $P$ is the 
unique parabolic subgroup of $G$ which contains $B$ such that the root space $\lieg_{-\alpha}$
for $\alpha \in S$ is contained in the Lie algebra of $P$ if and only if $\alpha \in S_P$.
The Lie algebra of $L^{(P)}$ is generated by the root spaces $\lieg_\alpha$ and $\lieg_{-\alpha}$
for $\alpha \in S_P$ together with the Lie algebra $\liet_c = \liet \otimes_\RR \CC$ of the complexification
$T_c$ of $T$. In addition the Lie algebra of $U$ is
\begin{equation} \lieu = \bigoplus_{\alpha \in R^+: \lieg_\alpha \not\subseteq {\rm Lie}(L^{(P)})} \lieg_\a \end{equation}
where $R^+$ is the set of positive roots for $G$, while the Lie algebra of $P$ is
\begin{equation} \liep = \liet_c \oplus \bigoplus_{\a \in R(S_P)} \lieg_\a \end{equation}
where $R(S_P)$ is the union of $R^+$ with the set of all roots which
can be written as sums of negatives of the simple roots in $S_P$.  If we identify
$S$ with the set of vertices of the Dynkin diagram of $K$ then the Dynkin diagram of the semisimple
part $Q^{(P)} = [K^{(P)},K^{(P)}]$ of $K^{(P)}$ is the subdiagram given by leaving out the vertices
which do not belong to $S_P$. We can decompose $\liek^{(P)} = {\rm Lie} K^{(P)}$ and 
$\liet$ as
$$\liek^{(P)} = [\liek^{(P)},\liek^{(P)}] \oplus \liez^{(P)} \ \ \mbox{ and } \ \ \liet = \liet^{(P)} \oplus
\liez^{(P)}$$
where $[\liek^{(P)},\liek^{(P)}]$ is the Lie algebra of $Q^{(P)} = [K^{(P)},K^{(P)}]$, while $\liet^{(P)}$ is the Lie
algebra of the maximal torus $T^{(P)} = T \cap [K^{(P)},K^{(P)}]$ of $Q^{(P)}$, and $\liez^{(P)}$ is
the Lie algebra of the centre $Z(K^{(P)})$ of $K^{(P)}$.
As before let $B^{{\rm op}} = T_c U^{{\rm op}}_{{\rm max}}$ be the Borel subgroup of $G$, with unipotent
radical $U_{{\rm max}}^{{\rm op}}$, which is opposite to $B$ in the sense that $B \cap B^{{\rm op}}
=T_c$ and $\umax \cap U_{{\rm max}}^{{\rm op}} = \{ 1 \}$, and let $\iota:\liets \to \liets$ be the involution given by
$\iota \l = - w_0 \l$ where $w_0$ denotes the element of the Weyl group $W$ of $G$ such that
$w_0 \umax w_0^{-1} = U_{{\rm max}}^{{ \rm op}}$.

By \cite{Grosshans2} Theorem 2.2 $U$ is a Grosshans subgroup of $G$, and so,
just as in the case when $U=\umax$, the ring of invariants $\calo(G)^U$ is 
finitely generated and $G/U$ has a canonical affine completion
\begin{equation} G/U \subseteq \overline{G/U}^{{\rm aff}} = {\rm Spec}(\calo(G)^U) \end{equation}
such that the complement of $G/U$ in $\overline{G/U}^{{\rm aff}}$
has codimension two. 

\begin{rem} \label{Iwas}
When $U = \umax$ the Iwasawa decomposition
$$G = K \ \exp(i\liet) \ \umax$$
enables us to identify $G/\umax$ with $K\exp(i\liet)$. More generally
we have an analogous decomposition
\begin{equation} G = K \times_{K^{(P)}} P = K \times_{K^{(P)}} \ L^{(P)} U = K \times_{K^{(P)}} \ K^{(P)} \exp(i \liek^{(P)}) U
= K \exp(i\liek^{(P)}) U \end{equation}
which enables us to identify $G/U$ with $K \exp(i \liek^{(P)})$.
\end{rem}

Let $X$ be a complex projective
variety on which $G$ acts linearly with respect to a very ample
line bundle $L$, and let $\omega$ be an associated $K$-invariant K\"{a}hler form on
$X$. Then it follows by the Borel transfer theorem (see e.g. \cite{Dolg} Lemma 4.1) that 
$$\hat{\calo}_L(X)^U \cong (\hat{\calo}_L(X) \otimes \calo(G)^U)^G$$
is finitely generated, and the associated projective variety
$$X/\!/U = {\rm Proj}(\hat{\calo}_L(X)^U)$$
is isomorphic to the GIT quotient $ (\overline{G/U}^{{\rm aff}} \times X)/\!/G$.
Just as in the case when $U=\umax$, if we have a suitable $K$-invariant K\"{a}hler
form on $\overline{G/U}^{{\rm aff}}$, then we will be able to identify $\xu$ with a symplectic quotient of 
$\overline{G/U}^{{\rm aff}} \times X$ by $K$, and obtain a symplectic description of $\xu$ analogous
to symplectic implosion, with $\overline{G/U}^{{\rm aff}}$ playing the r\^{o}le of the universal
imploded cross-section $\tkimp$. As is observed in \cite{GJS} $\S$6, the easiest case is when $K$
is semisimple and simply connected (for example when $K=SU(r+1)$); for general compact connected
$K$ one can reduce to this case by considering the product $\tilde{K}$ of the centre of $K$ and the universal cover of its commutator subgroup $[K,K]$, and expressing $K$ as $\tilde{K}/\Upsilon$, where $\Upsilon$ is a finite
central subgroup of $\tilde{K}$. 

Therefore, as in the previous subsection, let $K$ be a semisimple, connected and simply connected compact group,
let $\Lambda = \ker(\exp |_{\liet})$ be the exponential lattice in $\liet$, and let
$\Lambda^* = {\rm Hom}_{\ZZ}(\Lambda,\ZZ)$ be the weight lattice in $\liets$, so that 
$\Lambda^*_+ = \Lambda^* \cap \liets_+$ is the monoid of dominant weights. For $\lambda \in \Lambda^*_+$
let $V_{\lambda}$ be the irreducible $G$-module with highest weight $\lambda$, and let
$\Pi = \{\varpi_1, \ldots, \varpi_r \}$
be the set of fundamental weights, forming a $\ZZ$-basis of $\Lambda^*$ and a minimal set of generators
for $\Lambda^*_+$. 
Recall that 
we have an isomorphism of
$G \times G$-modules
\begin{equation} \label{thisiso3}
\calo(G) \cong \bigoplus_{\lambda \in \Lambda^*_+} V_{\lambda} \otimes V_\l^* 
\cong \bigoplus_{\lambda \in \Lambda^*_+} V_{\lambda} \otimes V_{\iota \l} \end{equation}
which restricts to  an isomorphism of $G \times T_c$-modules
\begin{equation} \label{thisiso2}
\calo(G)^\umax \cong \bigoplus_{\lambda \in \Lambda^*_+}  V_\l^{(T)} \otimes V_{\lambda}^* 
\cong \bigoplus_{\lambda \in \Lambda^*_+} V_{\lambda}^* \end{equation}
which  is generated as an algebra
 by its finite-dimensional vector subspace
$$E^* = \bigoplus_{\varpi \in \Pi} V_{\varpi}^*$$
giving us a closed
$G \times T_c$-equivariant embedding of $\overline{G/U}_{{\rm max}}^{{\rm aff}} = {\rm Spec}(\calo(G)^\umax)$
into the affine space $E$ equipped with
 a flat K\"{a}hler
structure. 
We have seen how Guillemin, Jeffrey and Sjamaar identify $\tkimp$ with 
$\overline{G/U}_{{\rm max}}^{{\rm aff}}$
 equipped with the K\"{a}hler structure obtained from this embedding in $E$.
To extend their construction to $\overline{G/U}^{{\rm aff}}$ when $U$ is the unipotent radical
of a parabolic subgroup $P \geq B$ as above, we first observe
from the proof of \cite{Grosshans2} Theorem 2.2 that $\calo(G)^U$ is generated by any finite-dimensional
${L^{(P)}}$-invariant (or equivalently ${K^{(P)}}$-invariant) vector subspace of 
$$\calo(G)  \cong \bigoplus_{\lambda \in \Lambda^*_+} V_{\lambda} \otimes V_{ \l}^*
\cong \bigoplus_{\lambda \in \Lambda^*_+} V_{\lambda} \otimes V_{\iota \l}
$$ which contains 
$$E^* = \bigoplus_{\varpi \in \Pi} V_{\varpi}^* \cong \bigoplus_{\varpi \in \Pi}
 V_{\varpi}^{(T)} \otimes V_{\varpi}^*.$$
Here as above $V_{ \varpi}^{(T)}$ is the irreducible $T_c$-module with weight $\varpi$
while
${K^{(P)}}=K \cap {L^{(P)}} = K \cap P$ is a maximal compact subgroup of the 
Levi subgroup ${L^{(P)}} = {K^{(P)}}_c$ of $P$, and ${K^{(P)}}$ acts on $\calo(G)$ via left multiplication on $G$.

Let $E^{{(P)}}$ be the dual of the smallest ${K^{(P)}}$-invariant subspace $(E^{{(P)}})^*$ 
of $\calo(G)$ containing $E^*$;
then $(E^{{(P)}})^*$ is fixed pointwise by $U$ since $K^{(P)}$ normalises $U$ and
$U$ is a subgroup of $\umax$ which fixes $E$ pointwise. The inclusion of $(E^{{(P)}})^*$ in $\calo(G)^U \subseteq \calo(G)$ induces a closed $ {L^{(P)}}
\times G$-equivariant 
embedding of $\overline{G/U}^{{\rm aff}} = {\rm Spec}(\calo(G)^U)$ into the affine space $E^{{(P)}}$,
whose projection to $E$ induces the embedding of $\overline{G/U}_{{\rm max}}^{{\rm aff}}$ described in the previous
subsection.
%

$({E^{(P)}})^*$ decomposes under the action of $K \times {K^{(P)}}$ as a direct sum of irreducible
$K \times {K^{(P)}}$-modules
$$({E^{(P)}})^* = \bigoplus_{\varpi \in \Pi} (V_{\varpi}^{({{P}})})^*$$
where $(V_{\varpi}^{({P})})^*$ is the smallest $K \times {K^{(P)}}$-invariant 
 subspace of $\calo(G)$ containing
$V_{\varpi}^*$. 
As in \cite{Grosshans} $\S$12 we have
$$(V_{\varpi}^{({{P}})})^* \cong  V_{\varpi}^{{K^{(P)}}} \otimes V_{\varpi}^* $$
where $V_{\varpi}^{{K^{(P)}}}$ is the irreducible $K^{(P)}$-module with highest
weight $\varpi$, 
so
\begin{equation} {E^{(P)}} = \bigoplus_{\varpi \in \Pi} V_{\varpi}^{({{P}})}
= \bigoplus_{\varpi \in \Pi}  (V_{ \varpi}^{K^{(P)}})^* \otimes
V_{\varpi} .\end{equation}
Moreover, if 
 $v^{(P)}_{\varpi}$ is the vector in $V_{ \varpi}^{({{P}})} \cong
(V_{\varpi}^{K^{(P)}})^* \otimes V_{\varpi}$ representing the inclusion of 
$V_{\varpi}^{K^{{(P)}}}$ in $V_{\varpi}$ then
the embedding  of $G/U \subseteq \overline{G/U}^{{\rm aff}}$ in $E^{(P)}$ induced by the inclusion of $(E^{(P)})^*$ in $\calo(G)^U$ 
takes the identity coset $U$ to $\sum_{\varpi \in \Pi} v_{\varpi}^{(P)}$.
Let
$$V_{\varpi}^{{K^{(P)}}} = \bigoplus_{\lambda \in \Lambda_{\varpi}^*} V_{\varpi,\lambda}^{{K^{(P)}}}$$
be the decomposition of $V_{\varpi}^{{K^{(P)}}}$ into weight spaces with weights
$\lambda \in \liets$ under the action of the maximal torus $T$ of $K^{(P)}$. Then 
$V_{\varpi}^{{{(P)}}}$ decomposes as a $K \times T$-module into a sum of 
irreducible $K \times T$-modules
\begin{equation} V_{\varpi}^{{{(P)}}} \cong \bigoplus_{\lambda} V_{\varpi} \otimes (V_{\varpi,\lambda}^{{K^{(P)}}})^* \end{equation}
and  
$v_{\varpi}^{(P)} = \sum_{\lambda} v_{\varpi,\lambda}^{(P)}$
where $v_{\varpi,\lambda}^{(P)} \in V_{\varpi} \otimes (V_{\varpi,\lambda}^{{K^{(P)}}})^*$
represents the inclusion of $V_{\varpi,\lambda}^{{K^{(P)}}}$ in $V_{\varpi}$. In particular
$v_{\varpi,\varpi}^{(P)}$ is a highest  weight vector for the action of
$K \times K^{(P)}$ on $V^{(P)}_{\varpi}$. 

\begin{rem} \label{onto} The embedding  of $G/U \subseteq \overline{G/U}^{{\rm aff}}$ in $E^{(P)}$ induced by the inclusion of $(E^{(P)})^*$ in $\calo(G)^U$ 
takes the identity coset to $\sum_{\varpi \in \Pi} v_{\varpi}^{(P)}$. From the 
decomposition $G = K \exp(i\liek^{(P)}) U$ (see Remark \ref{Iwas} above) and the
compactness of $K$ it follows that the closure $\overline{G/U}^{{\rm aff}}$ of the
$G$-orbit of $\sum_{\varpi \in \Pi} v_{\varpi}^{(P)}$ in $E^{(P)}$ is given by the
$K$-sweep
$$\overline{G/U}^{{\rm aff}} = K (\overline{\exp (i \liek^{(P)}) \sum_{\varpi \in \Pi} v^{(P)}_{\varpi}})$$
of the closure in $E^{(P)}$ of the $\exp (i \liek^{(P)})$-orbit of $\sum_{\varpi \in \Pi} v^{(P)}_{\varpi}$.
Similarly the closure in $E^{(P)}$ (or equivalently in the
linear subspace $\bigoplus_{\varpi \in \Pi}  (V_{ \varpi}^{K^{(P)}})^* \otimes
V_{\varpi}^{K^{(P)}} $ of $E^{(P)}$)
of the $L^{(P)}$-orbit of $\sum_{\varpi \in \Pi} v^{(P)}_{\varpi}$
(which is a free orbit since $U \cap L^{(P)} = \{ 1\}$) is given by
$K^{(P)} (\overline{\exp (i \liek^{(P)}) \sum_{\varpi \in \Pi} v^{(P)}_{\varpi}})$.
Note alse that
$\liek^{(P)} = \bigcup_{k \in K^{(P)}} {\rm Ad}(k)\liet$
and so
\begin{equation} \exp (i \liek^{(P)}) = \bigcup_{k \in K^{(P)}} k \ \exp(i\liet) \ k^{-1}. \end{equation}
\end{rem}

Let $S_{P} = \{\a_1,\ldots \a_{r(P)}\} \subseteq S = \{\a_1,\ldots \a_{r}\} $ be the
set of simple roots for the root system of $(K^{(P)},T)$ with corresponding positive Weyl chamber
$$\liets_{+,P} = \{ \zeta \in \liets: \zeta \cdot \a \geq 0 \mbox{ for all }\a \in S_{P} \}
= \liet^{(P)*}_+ \oplus \liez^{(P)*}$$
where $\liez^{(P)}$ is the Lie algebra of the centre $Z(K^{(P)})\leq T $ of $K^{(P)}$ and
$\liet^{(P)*}_+$ is the positive Weyl chamber for the semisimple part
$$ Q^{(P)} = [{K^{(P)}},K^{(P)}]$$ of $K^{(P)}$ with respect to the maximal torus 
$T^{(P)} = T \cap [K^{(P)},K^{(P)}]$ of $Q^{(P)}$ and simple roots given by restricting $S_{P}$ to $T^{(P)}$.
If $\varpi \cdot \alpha = 0$ for all $\alpha \in S_{P}$ (or equivalently
if $\varpi = \varpi_j$ for $j > r(P)$) then $\varpi \in \liez^{(P)*}$ and 
$V_{\varpi}^{{K^{(P)}}}$ is one-dimensional; in this situation $Q^{(P)}$ acts
trivially on $V_{\varpi}^{{K^{(P)}}}$ and we have $V_{\varpi}^{{K^{(P)}}}
= V_{\varpi,\varpi}^{{K^{(P)}}}$ with $V_{\varpi,\lambda}^{{K^{(P)}}}=0$ if
$\lambda \neq \varpi$, and $v_{\varpi}^{{K^{(P)}}}
= v_{\varpi,\varpi}^{{K^{(P)}}}$ while $v_{\varpi,\lambda}^{{K^{(P)}}}=0$ if
$\lambda \neq \varpi$. On the other hand if $j \leq r(P)$ then $\varpi = \varpi_j$ restricts to
a fundamental weight for $Q^{(P)}$ and $V_{\varpi}^{{K^{(P)}}} = V_{\varpi}^{{Q^{(P)}}}$
is the irreducible $Q^{(P)}$-module with highest weight $\varpi|_{Q^{(P)}}$ on
which $Z(K^{(P)})$ acts
as scalar multiplication by $\varpi|_{Z(K^{(P)})}$.

There is a unique $K \times {K^{(P)}}$-invariant Hermitian inner product on ${E^{(P)}} = \bigoplus_{\varpi \in \Pi}
V_{\varpi}^{(P)}$
satisfying $|\!| v_{\varpi,\varpi}^{(P)}|\!| = 1$ for each $\varpi \in \Pi$, which is obtained from
$K$-invariant Hermitian inner products on the irreducible $K$-modules $V_{\varpi}$
and their restrictions to 
$K^{(P)}$-invariant Hermitian inner products on the irreducible $K^{(P)}$-modules $V_{\varpi}^{(P)}$.
This gives ${E^{(P)}}$ a flat
K\"{a}hler structure which is $K \times {K^{(P)}}$-invariant. 

\begin{rem} \label{star}
Recall that $${E^{(P)}} = \bigoplus_{\varpi \in \Pi} V_{\varpi}^{({{P}})}
= \bigoplus_{\varpi \in \Pi}  (V_{ \varpi}^{K^{(P)}})^* \otimes
V_{\varpi} $$
where $V_{\varpi}^{K^{{(P)}}} \subseteq V_{\varpi}$, and 
the embedding  of $G/U \subseteq \overline{G/U}^{{\rm aff}}$ in $E^{(P)}$ induced by the inclusion of $(E^{(P)})^*$ in $\calo(G)^U$ 
takes the identity coset $U$ to $\sum_{\varpi \in \Pi} v_{\varpi}^{(P)}$ where
 $v^{(P)}_{\varpi} \in V_{ \varpi}^{({{P}})} \cong
(V_{\varpi}^{K^{(P)}})^* \otimes V_{\varpi}$ represents the inclusion of 
$V_{\varpi}^{K^{{(P)}}}$ in $V_{\varpi}$. Thus
$$\sum_{\varpi \in \Pi} v_{\varpi}^{(P)} \in \bigoplus_{\varpi \in \Pi}  (V_{ \varpi}^{K^{(P)}})^* \otimes
V_{\varpi}^{K^{(P)}} \subseteq {E^{(P)}}$$
where $\bigoplus_{\varpi \in \Pi}  (V_{ \varpi}^{K^{(P)}})^* \otimes
V_{\varpi}^{K^{(P)}}$ is invariant under the action of the subgroup $K^{(P)} \times K^{(P)}$ of
$K \times K^{(P)}$ on $E^{(P)}$, and indeed is invariant under the action of $L^{(P)} \times L^{(P)}$.
If we identify $(V_{ \varpi}^{K^{(P)}})^* \otimes
V_{\varpi}^{K^{(P)}}$ with ${\rm End}(V_{ \varpi}^{K^{(P)}})$ equipped with the
Hermitian structure 
$$\langle A,B \rangle = {\rm Trace}(AB^*)$$
in the standard way, then $v^{(P)}_{\varpi}$ is identified with the identity map
in ${\rm End}(V_{ \varpi}^{K^{(P)}})$. If $V$ is any Hermitian vector space then the moment map
for the action of the product of unitary groups 
$U(V) \times U(V)$  on ${\rm End}(V)$ 
by left and right multiplication 
is given (up to a nonzero real
 scalar) by
 $$A \mapsto (iAA^*,iA^*A)$$
 (cf. \cite{Paradan} $\S$3.3). Thus the moment map for the action of $K^{(P)} \times K^{(P)}$ on
$$ \bigoplus_{\varpi \in \Pi}  (V_{ \varpi}^{K^{(P)}})^* \otimes
V_{\varpi}^{K^{(P)}} \cong \bigoplus_{\varpi \in \Pi}  {\rm End}(V_{ \varpi}^{K^{(P)}})$$
 is given (up to multiplication by a nonzero real
 scalar) by
 \begin{equation} \sum_{\varpi \in \Pi} A_{\varpi} \mapsto (\pi^{K^{(P)}}(
 \sum_{\varpi \in \Pi} iA_{\varpi}A_{\varpi}^*),\pi^{K^{(P)}}( 
 \sum_{\varpi \in \Pi} iA_{\varpi}^*A_{\varpi})) \end{equation}
 where $\pi^{K^{(P)}}:\lieu(\bigoplus_{\varpi \in \Pi}  
V_{\varpi}^{K^{(P)}})^* \to (\liek^{(P)})^*$ is the projection induced
by the inclusion of $K^{(P)}$ as a subgroup of the
unitary group $U(\bigoplus_{\varpi \in \Pi}  
V_{\varpi}^{K^{(P)}})$. In particular if $g$ belongs to the complexification
$L^{(P)}$ of $K^{(P)}$ and $g_{\varpi}:V_{\varpi}^{K^{(P)}} \to V_{\varpi}^{K^{(P)}}$
is the action of $g$ on $V_{\varpi}^{K^{(P)}}$, then
$$ g \sum_{\varpi \in \Pi} v_{\varpi}^{(P)} = \sum_{\varpi \in \Pi} g_{\varpi}$$
and the moment map for the left $K^{(P)}$-action sends this to
$$\pi^{K^{(P)}}(
 \sum_{\varpi \in \Pi} ig_{\varpi}g_{\varpi}^*) \in  \liek^{(P)}.$$
Using the decomposition $\liek^{(P)} = [\liek^{(P)},\liek^{(P)}] \oplus \liez^{(P)}$ we can decompose
 $\pi^{K^{(P)}}:\lieu(\bigoplus_{\varpi \in \Pi}  
V_{\varpi}^{K^{(P)}})^* \to (\liek^{(P)})^*$ as
\begin{equation} \pi^{K^{(P)}}= \pi^{[K^{(P)},K^{(P)}]} \oplus \pi^{Z(K^{(P)})}:\lieu(\bigoplus_{\varpi \in \Pi}  
V_{\varpi}^{K^{(P)}})^* \to [\liek^{(P)},\liek^{(P)}]^* \oplus \liez^{(P)*}.
\end{equation}
If $g=yz$ with $y \in [L^{(P)},L^{(P)}] = Q^{(P)}_c$ and $z \in Z(L^{(P)})=
Z(K^{(P)})_c$, then the $K^{(P)}$-moment map above sends 
$ g \sum_{\varpi \in \Pi} v_{\varpi}^{(P)}$ to 
$$\pi^{[K^{(P)},K^{(P)}]}(
 \sum_{1\leq j \leq r(P)} iy_{\varpi_j}y_{\varpi_j}^*)
 + \pi^{Z(K^{(P)})}(\sum_{1\leq j \leq r} iz_{\varpi_j}z_{\varpi_j}^*)
 \in  [\liek^{(P)},\liek^{(P)}]^* \oplus \liez^{(P)*}
 .$$
 It follows by the arguments of \cite{Paradan} $\S$3 (in particular Proposition
 3.10) that the $T^{(P)}_c$-orbit of $ \sum_{\varpi \in \Pi} v_{\varpi}^{(P)}$
 is mapped diffeomorphically onto $\liet^{(P)}$ by the moment map
 \begin{equation} y \sum_{\varpi \in \Pi} v_{\varpi}^{(P)} \mapsto \pi^{T^{(P)}}(\sum_{1\leq j \leq r(P)} iy_{\varpi_j}y_{\varpi_j}^*) \end{equation}
for the action of $T^{(P)}$ on $E^{(P)}$, since its image in the projective
space $\PP(E^{(P)})$ is mapped diffeomorphically by the associated moment map
onto the
convex hull of the set $\{ w\varpi: \varpi \in \Pi, w \in W^{(P)} \}$ where $W^{(P)}$
is the Weyl group of $Q^{(P)} = [K^{(P)},K^{(P)}]$ (cf. Remark \ref{calffo}).
\end{rem}

Now consider the moment map $\mu^{E^{(P)}}_T$ for the restriction to $T$ of the $K^{(P)}$-action
on $E^{(P)}$. This is given (up to multiplication by a positive constant) by
$$\sum_{\varpi,\lambda} u_{\varpi,\lambda} \mapsto \sum_{\varpi,\lambda} |\!| u_{\varpi,\lambda}|\!|^2 \lambda$$
when $u_{\varpi,\lambda} \in V_{\varpi,\lambda}^{{K^{(P)}}}$ for $\varpi \in \Pi $ and $\lambda \in \Lambda^*_\varpi \subseteq \Lambda^*$.
The embedding  of $G/U \subseteq \overline{G/U}^{{\rm aff}}$ in $E^{(P)}$ induced by the inclusion of $(E^{(P)})^*$ in $\calo(G)^U$ 
takes  the coset of $t\in T_c$ to
$$\sum_{\varpi,\lambda} \lambda(t)^{-1} v_{\varpi,\lambda}^{(P)},$$
so the value taken by this moment map on the coset $tU$ of $t = t_1 t_2 \in T_c$ 
where $t_1 \in T^{(P)}_c$ and $t_2 \in Z(L^{(P)}) = Z(K^{(P)})_c$ is given by
\begin{equation} \label{sr}
\sum_{\varpi,\lambda} |\lambda(t)|^{-2} |\!| v^{(P)}_{\varpi,\lambda}|\!|^2 \lambda=
\sum_{j=1}^{r(P)}|\varpi_j(t_2)|^{-2} \sum_{\lambda} |\lambda(t_1)|^{-2} |\!| v^{(P)}_{\varpi_j,\lambda}|\!|^2 \lambda
+ \sum_{j=r(P)+1}^r |\varpi_j(t_2)|^{-2} |\!| v^{(P)}_{\varpi_j,\varpi_j}|\!|^2 \varpi_j
\end{equation} 
where the $j$th sum over $\lambda$ runs over all the weights of the irreducible
$K^{(P)}$-module $V_{\varpi_j}^{K^{(P)}}$ with highest weight $\varpi_j$.
When we decompose $\liets$ as $\liet^{(P)*} \oplus \liez^{(P)*}$ this has
component $\sum_{j=1}^r |\varpi_j(t_2)|^{-2} |\!| v^{(P)}_{\varpi_j}|\!|^2 \varpi_j|_{Z(K^{(P)}}$
in $\liez^{(P)*}$ and $\sum_{j=1}^{r(P)}|\varpi_j(t_2)|^{-2} \sum_{\lambda} |\lambda(t_1)|^{-2} |\!| v^{(P)}_{\varpi_j,\lambda}|\!|^2 \lambda|_{T^{(P)}}$ in $\liet^{(P)*}$.

\begin{defn} Let $\tcone $ be the cone
$$\tcone  
= \bigcup_{w \in W^{(P)}} {\rm Ad}^*(w) \liets_+$$
in $\liets$, 
where $W^{(P)}$ is the Weyl group of $Q^{(P)}=[K^{(P)},K^{(P)}]$ (which is a subgroup of
the Weyl group $W$ of $K$). 
\end{defn}


\begin{lemma} \label{whh}
The restriction to the closure $\overline{\exp(i\liet) \sum_{\varpi \in \Pi} v^{(P)}_{\varpi}}$
of the $\exp(i\liet)$-orbit in $E^{(P)}$ of $ \sum_{\varpi \in \Pi} v^{(P)}_{\varpi}$
of  the moment map $\mu^{E^{(P)}}_T$ for the action of $T$ 
on $E^{(P)}$ is a homeomorphism onto the cone
$\tcone $
in $\liets$. Its inverse 
provides a continuous injection
\begin{equation} \label{calffP} \calf^{(P)} : \tcone  \to \overline{G/U}^{{\rm aff}} \subseteq E^{(P)} \end{equation}
such that $\mu_T^{E^{(P)}} \circ \calf^{(P)}$ is the identity on $\tcone $.
Moreover $\overline{\exp(i\liet) \sum_{\varpi \in \Pi} v_\varpi^{(P)}}$
is the union of finitely many $\exp(i\liet)$-orbits, each of the form $$
\calf^{(P)}(\s) = \exp(i\liet) \sum_{\varpi \in \Pi, \lambda \in \Lambda^*_\varpi \cap \bar{\s}} v^{(P)}_{\varpi,\lambda}$$
where $\s$ 
is an open face of $\tcone $.
\end{lemma}
\noindent{\bf Proof.} 
This follows by applying the results of
\cite{Atiyah} to the compactification $\PP(\CC \oplus E^{(P)})$ of the affine
space $E^{(P)}$, as in Remark \ref{calffo},
and observing that the convex hull of the weights $\lambda$ of the $T$-action on the $K^{(P)}$-module $V^{K^{(P)}}_\varpi$ is the convex hull of 
$\{ w\varpi: w \in W^{(P)}\}$, and thus the convex hull of the half-lines
$\RR_+\lambda$ for $\lambda \in \Lambda_\varpi$ with $\varpi \in \Pi$ is the
cone $\tcone $. \hfill $\Box$

\begin{lemma} (cf. \cite{Paradan} Lemma 3.12) \label{Parad}
The image of the closure $\overline{T_c \sum_{\varpi \in \Pi} v^{(P)}_{\varpi}}$
of the $T_c$-orbit in $E^{(P)}$ of $ \sum_{\varpi \in \Pi} v^{(P)}_{\varpi}$ under the
$K^{(P)}$-moment map $\mu^{E^{(P)}}:E^{(P)} \to (\liek^{(P)})^* \cong \liek^{(P)}$ is
contained in $\liet$.
\end{lemma}
\noindent{\bf Proof:} The orthogonal complement to $\liet$ in $\liek^{(P)}$ is
$[\liek^{(P)},\liet]$, and if $\zeta \in \liet$ and $\xi \in \liek^{(P)}$ and $t \in T_c$
then by Remark \ref{star}
$$\mu^{E^{(P)}}(t  \sum_{\varpi \in \Pi} v^{(P)}_{\varpi}).[\xi,\zeta]
= \sum_{\varpi \in \Pi} {\rm Trace}(i[\xi,\zeta]t_{\varpi} t_{\varpi}^*)
= \sum_{\varpi \in \Pi} {\rm Trace}(i\xi[\zeta,t_{\varpi} t_{\varpi}^*]) =0$$
since $[\zeta,t_{\varpi} t_{\varpi}^*]=0.$  \hfill $\Box$

\begin{corollary} \label{thecor}
The restriction of  the
$K^{(P)}$-moment map $\mu^{E^{(P)}}:E^{(P)} \to (\liek^{(P)})^*$ to the closure
$$\overline{\exp(i\liek^{(P)}) \sum_{\varpi \in \Pi} v^{(P)}_{\varpi}}$$
of the $\exp(i\liek^{(P)})$-orbit in $E^{(P)}$ of $ \sum_{\varpi \in \Pi} v^{(P)}_{\varpi}$
is a homeomorphism from $\overline{\exp(i\liek^{(P)}) \sum_{\varpi \in \Pi} v^{(P)}_{\varpi}}$
onto the closed subset
$$\liek^{(P)*}_+ =  {\rm Ad}^*(K^{(P)})\  \tcone $$
of $\liek^{(P)*}$.
Moreover $\overline{\exp(i\liek^{(P)}) \sum_{\varpi \in \Pi} v_\varpi^{(P)}}$
is the union of finitely many $\exp(i\liek^{(P)})$-orbits which
correspond under this homeomorphism to the open faces of $\liek^{(P)*}_{+}$.
\end{corollary}
\noindent{\bf Proof:} We have already observed that the restriction of the $T$-moment map 
$\mu_T^{E^{(P)}}:E^{(P)} \to \liets$ to the closure
$$\overline{\exp (i \liet) \sum_{\varpi \in \Pi} v^{(P)}_{\varpi}}$$
of the $\exp (i\liet)$-orbit of the image $\sum_{\varpi \in \Pi} v^{(P)}_{\varpi}$
of the identity coset $U$ under the embedding of $G/U$ in $E^{(P)}$ is a homeomorphism
from this closure onto the cone $\tcone $. Since $\mu_T^{E^{(P)}}$ is the projection
of $\mu^{E^{(P)}}$ onto $\liets$, it follows immediately from Lemma \ref{Parad} above that
the restriction of $\mu^{E^{(P)}}:E^{(P)} \to (\liek^{(P)})^* \cong \liek^{(P)}$
to this closure
$\overline{\exp (i \liet) \sum_{\varpi \in \Pi} v^{(P)}_{\varpi}}$
 is a homeomorphism onto the cone $\tcone $ when $\liets$
 is identified with $\liet \subseteq \liek^{(P)}$ via the restriction of the fixed
 invariant inner  product on $\liek$. Replacing the maximal torus $T$ with
 $kTk^{-1}$ for any $k \in K^{(P)}$ it follows that
 the restriction of $\mu^{E^{(P)}}:E^{(P)} \to (\liek^{(P)})^* $
to the closure
$\overline{k\exp (i \liet)k^{-1} \sum_{\varpi \in \Pi} v^{(P)}_{\varpi}}$
of the $\exp (i{\rm Ad}(k)\liet)$-orbit of the image $\sum_{\varpi \in \Pi} v^{(P)}_{\varpi}$
of the identity coset $U$ under the embedding of $G/U$ in $E^{(P)}$
 is a homeomorphism onto the cone ${\rm Ad}^*(k)\tcone $. 
  Putting these homeomorphisms together for $k \in K^{(P)}$ we get a homeomorphism
$\calm$ from 
$$\calz = \{ (kN^{(P)}_T, x) \in K^{(P)}/N^{(P)}_T \ \times \ E^{(P)}\ \  :\ \  x \in
\overline{k\exp (i \liet)k^{-1} \sum_{\varpi \in \Pi} v^{(P)}_{\varpi}} \},$$
where $N^{(P)}_T$ is the normaliser of $T$ in $K^{(P)}$, to 
$$K^{(P)} \times_{N^{(P)}_T} \tcone 
= \{ (k N^{(P)}_T,\xi) \in K^{(P)}/N^{(P)}_T \times \liek^{(P)*}\ \ |\ \  \xi \in {\rm Ad}^*(k) \tcone  \}$$
which fits into a diagram 
\begin{equation} \label{fpicture}
\begin{array}{ccccc}
   \calz & \rightarrow & K^{(P)} \times_{N^{(P)}_T} \tcone  & & \\
      &  &  &  & \\
  \alpha \Big\downarrow & & \Big\downarrow \beta & & \\
  &  &  &  &  \\
 \overline{\exp(i\liek^{(P)}) \sum_{\varpi \in \Pi} v^{(P)}_{\varpi}} & \rightarrow  & 
 \bigcup_{k \in K^{(P)}} {\rm Ad}^*(k)\  \tcone 
& = &  \liek^{(P)*}_+
     \end{array} \end{equation}
where the first horizontal map is the homeomorphism $\calm$ and  the second
 is $\mu^{E^{(P)}}$. Since the image of $\alpha$ is dense and $K^{(P)}$
is compact, it follows that $\alpha$ is surjective. Moreover $\beta$ is surjective, and
$\beta(k N^{(P)}_T,\xi)=\beta(k' N^{(P)}_T,\xi') $ if and only if 
${\rm Ad}^*(k^{-1})\xi$ lies in an open face $\sigma$ of $\liet_+$ such
that $k'k^{-1} \in K^{(P)}_\sigma$, in which case
$\alpha(\calm^{-1}(k N^{(P)}_T,\xi)) = \alpha(\calm^{-1}(k N^{(P)}_T,\xi))$.
Thus $$\mu^{E^{(P)}}:\overline{\exp(i\liek^{(P)}) \sum_{\varpi \in \Pi} v^{(P)}_{\varpi}} \to \liek^{(P)*}_+$$
is a continuous bijection, which is a homeomorphism since $K$ is compact and $\calm$ is a homeomorphism. \hfill $\Box$


The inverse of $\mu^{E^{(P)}}:\overline{\exp(i\liek^{(P)}) \sum_{\varpi \in \Pi} v^{(P)}_{\varpi}} \to \liek^{(P)*}_+$ gives us a continuous
$K^{(P)}$-equivariant map
$$ \calf^{(P)} : \liek_{+}^{(P)*} \to \overline{G/U}^{{\rm aff}} \subseteq E^{(P)} $$ extending (\ref{calffP})
such that $\mu_T^{E^{(P)}} \circ \calf^{(P)}$ is the identity on $\liek^{(P)*}_{+}$. This in turn extends to  a continuous
$K \times K^{(P)}$-equivariant map
\begin{equation} \label{calffP2} \calf^{(P)} : K \times \liek_{+}^{(P)*} \to \overline{G/U}^{{\rm aff}}  \end{equation}
which is surjective since $\overline{G/U}^{{\rm aff}} = K (\overline{\exp (i \liek^{(P)}) \sum_{\varpi \in \Pi} v^{(P)}_{\varpi}})$ by Remark \ref{onto}.




\begin{defn} \label{dffn}
If $ \zeta \in \liek^{(P)*}_+ = {\rm Ad}^*({K^{(P)}})\tcone = {\rm Ad}^*({K^{(P)}})\liets_{+}$
let $\zeta = {\rm Ad}^*(k)\xi$ with $k \in K^{(P)}$ and $\xi \in \liets_+$, and let $\s_0$ be the open face of 
$\liets_+$ containing $\xi$. Let $\s_0(P)$ be the open face of $\liets_+$ whose closure is
$$\overline{\sigma_0(P)}
= \{ \zeta \in \liets: \zeta \cdot \a = 0 \mbox{ for all }\a \in R_{\s_0} \setminus R^{(P)} \}$$
where $R$ and $R^{(P)}$ are the sets of roots of $K$ and $K^{(P)}$, and
$$ R_{\s_0} 
= \{ \alpha \in R: \zeta \cdot \a = 0 \mbox{ for all }\zeta \in \s_0 \},$$ 
so that $\s_0(P)$ is an open subset of the open face containing $\s_0$ of the
cone $\tcone$.
Finally let $K_\zeta(P) = k K_\xi k^{-1}$ where $K_\xi(P) = K_{\s_0(P)}$ is the stabiliser
under the adjoint action of $K$ of any element of $\s_0(P)$.
\end{defn}

Note that $K_\zeta(P) \leq K_\zeta$ for any $\zeta \in \liek^{(P)*}_+$. 

\begin{lemma} \label{GJS6.2}
(cf. \cite{GJS} Lemma
6.2)

Let $\sigma$ be an open face of $\tcone $ and let
$$v_\sigma^{(P)} = \sum_{\varpi \in \Pi, \lambda \in {\rm Ad}^*(W^{(P)})\varpi \cap \overline{\s}} v_{\varpi,\lambda}^{(P)}.$$
If $\zeta \in \sigma$ then the stabiliser of $v^{(P)}_{\s}$ in $K$ is
$[K_\zeta(P),K_\zeta(P)]$.
\end{lemma}
\noindent {\bf Proof:} 
Recall that $\tcone  = \bigcup_{w \in W^{(P)}} {\rm Ad}^*(w) \liets_+$, so there is an element $w_0$ of the Weyl group 
$W^{(P)}$ of $Q^{(P)} = [K^{(P)},K^{(P)}]$ such that ${\rm Ad}^*(w_0)\zeta \in \liets_+$ and
${\rm Ad}^*(w_0)\sigma$ contains an open face
$\s_0$ of $\liets_+$ with $\s_0$ an open subset of $\s$. First assume that $\xi = {\rm Ad}^*(w_0)\zeta$
lies in $\s_0$.
Then if $\varpi \in \Pi$ and $w \in W^{(P)}$ we have ${\rm Ad}^*(w) \varpi \in \overline{\sigma}$ if and
only if ${\rm Ad}^*(w) \varpi$ lies in the linear subspace of $\liets$ spanned by $\Pi \cap \overline{\s} = \Pi \cap \overline{\s_0}$, and
since $\xi \in \s_0$ this happens if and only if $K_\x \leq w K_\varpi w^{-1}$, so that
$$\bigcap_{\varpi \in \Pi, \lambda = w\varpi \in W^{(P)}\varpi \cap \overline{\s}} w K_\varpi w^{-1} = K_\xi.$$
As in the proof of \cite{GJS} Lemma 6.2 we find that  if $w \in W^{(P)}$ the 
stabiliser in $G=K_c$ of $[v_{\varpi,w\varpi}^{(P)}] \in \PP((V_\varpi^{K^{(P)}})^* \otimes V_\varpi)$ is $ wP_\varpi w^{-1}$, where $P_\varpi$ is the parabolic subgroup of $G$ associated to $\varpi$, and thus the 
stabiliser in $K$ of $v_\sigma^{(P)}$ is the conjugate by $w_0$ of
$$\{ k \in \bigcap_{\varpi \in \Pi, \lambda = w\varpi \in W^{(P)}\varpi \cap \overline{\s}} w K_\varpi w^{-1} \ \ : \ \  \tilde{\lambda}(g) = 1 \mbox{ for all } \tilde{\lambda} \in \Lambda^* \cap \bar{\s} \} 
$$
\begin{equation} = \{ k \in K_\xi \ \ : \ \  \tilde{\lambda}(g) = 1 \mbox{ for all } \tilde{\lambda} \in \Lambda^* \cap \bar{\s} \} =
[K_\xi,K_\xi] = [K_{\s_0},K_{\s_0}]. \end{equation}

In general if $\xi = {\rm Ad}^*(w_0)\zeta$
lies in $\liets_+ \cap {\rm Ad}^*(w_0)\sigma $ then there is a unique open face $\s_0$ of $\liets_+$
containing $\xi$. Let $\s_0(P)$ be as in Definition \ref{dffn}; then $\overline{\s} \cap \liets_+ =
\overline{\s_0(P)}$, and so by the previous paragraph the stabiliser of $v^{(P)}_{\s}$ in $K$ is
$$w_0 [K_{\s_0(P)},K_{\s_0(P)}] w_0^{-1}  = [K_\zeta(P),K_\zeta(P)].$$
\hfill $\Box$


Thus we extend the definition of the imploded cross-section $\ximp$ to
a ${K^{(P)}}$-imploded cross-section $\ximpq$ as follows.

\begin{defn} \label{impq} 
Let $(X,\omega)$ be a symplectic manifold on which $K$
acts with a moment map $\mu:X \to \lieks$.
As before let
\begin{equation} \liek^{(P)*}_+ = {\rm Ad}^*({K^{(P)}})\tcone = {\rm Ad}^*({K^{(P)}})\liets_{+}
  = {\rm Ad}^*({Q^{(P)}})\liets_{+}
\subseteq {\liek^{(P)*}}  \label{lieqsplus} \end{equation}
be the sweep of $\liets_{+}$ under the co-adjoint action of ${K^{(P)}}$ on $\lieks$, and let $\Sigma^{(P)}$ be the set of open faces of $\liek^{(P)*}_{+}$.
If $\zeta \in \liek^{(P)*}$ let $K_\zeta(P) $ be defined as in Definition \ref{dffn}.
The {\em ${K^{(P)}}$-imploded cross-section}
of $X$ is 
$$ \ximpq = \mu^{-1}(\liek^{(P)*}_+)/\approx_{K^{(P)}} $$
where  $x \approx_{K^{(P)}} y$ if and only if 
$\mu(x) = \mu(y) = \zeta \in 
\liek^{(P)*}_+$
 and
$x = \kappa y$ for some $\kappa \in 
[K_\zeta(P),K_\zeta(P)]$.

The {\em universal ${K^{(P)}}$-imploded cross-section} is the ${K^{(P)}}$-imploded cross-section
$$ \tkimpq = K \times \liek^{(P)*}_+ / \approx_{K^{(P)}} $$
for the cotangent bundle $T^*K \cong K \times \lieks$ with respect to the $K$-action induced from the 
right action of $K$ on itself. 

\end{defn}

\begin{theorem} \label{maint} The map 
$\calf^{(P)} : K \times \liek_{+}^{(P)*} \to \overline{G/U}^{{\rm aff}}  $ 
of (\ref{calffP2}) induces a 
$K \times K^{(P)}$-equivariant homeomorphism 
$$ \tkimpq = K \times \liek^{(P)*}_+ / \approx_{K^{(P)}} 
 \to
\overline{G/U}^{{\rm aff}} \subseteq {E^{(P)}}. $$
 Moreover under this
identification of $ K \times \liek^{(P)*}_+ / \approx_{K^{(P)}} $ with 
$\overline{G/U}^{{\rm aff}} \subseteq {E^{(P)}}$, the moment map
for the action of $K \times K^{(P)}$ on $E^{(P)}$ is induced by
the map $
(K  \times \liek^{(P)*}_+)/\approx_{K^{(P)}} \to
\lieks \times \liek^{(P)*}$
given by 
$$(k,\zeta) \mapsto (Ad^*(k)(\zeta),\zeta)).$$
\end{theorem}
 
\noindent{\bf Proof:} By Lemma \ref{GJS6.2} 
$\calf^{(P)}$ induces a 
continuous map
$ \tkimpq  \to
\overline{G/U}^{{\rm aff}} \subseteq {E^{(P)}}, $
which is surjective since $\overline{G/U}^{{\rm aff}} = K (\overline{\exp (i \liek^{(P)}) \sum_{\varpi \in \Pi} v^{(P)}_{\varpi}})$ by Remark \ref{onto}.
The map $
(K  \times \liek^{(P)*}_+)/\approx_{K^{(P)}} \to
\lieks \times \liek^{(P)*}$
given by 
$$(k,\zeta) \mapsto (Ad^*(k)(\zeta),\zeta))$$
is the composition of $\calf^{(P)} : K \times \liek_{+}^{(P)*} \to \overline{G/U}^{{\rm aff}}  $ with the restriction to 
$\overline{G/U}^{{\rm aff}} \subseteq {E^{(P)}}$ of  the moment map
$\mu^{E^{(P)}}$ for the action of $K \times K^{(P)}$ on $E^{(P)}$.
Moreover  $\calf^{(P)}$ is continuous and surjective and
restricts to a homeomorphism from $\overline{\exp(i\liek^{(P)}) \sum_{\varpi \in \Pi} v^{(P)}_{\varpi}}$
to $\liek^{(P)*}_+$ by Corollary \ref{thecor}. If $\calf^{(P)}(k_1,\zeta_1) =
\calf^{(P)}(k_2,\zeta_2)$ then it follows by applying $\mu^{E^{(P)}}$ that
$(Ad^*(k_1)(\zeta_1),\zeta_1) = (Ad^*(k_2)(\zeta_2),\zeta_2)$ and therefore
$\zeta_1 = \zeta_2$ and $k_1 k_2^{-1} \in K_{\zeta_1} = K_{\zeta_2}$. Thus
$$\calf^{(P)}(1,\zeta_1) =(k_1^{-1},1)\calf^{(P)}(k_1,\zeta_1)
=(k_1^{-1},1) \calf^{(P)}(k_2,\zeta_2) = (k_1^{-1} k_2,1) \calf^{(P)}(1,\zeta_2) 
 = (k_1^{-1} k_2,1) \calf^{(P)}(1,\zeta_1).$$
Since $\zeta_1 = \zeta_2 \in \liek^{(P)*}_+ = {\rm Ad}^*(K^{(P)})\tcone $
we can write $\zeta_1 = \zeta_2 = {\rm Ad}^*(k_0)(\zeta)$ where $\zeta \in \tcone $ and $k_0 \in K^{(P)}$, so 
$$\calf^{(P)}(1,\zeta) =(1,k_0^{-1})\calf^{(P)}(1,\zeta_1)
=(k_1^{-1}k_2,k_0^{-1}) \calf^{(P)}(1,\zeta_1) = (k_1^{-1} k_2,1) \calf^{(P)}(1,\zeta).$$
By Lemma \ref{whh} $\calf^{(P)}(1,\zeta)$ lies in the
$\exp(i\liet)$-orbit of $\sum_{\varpi \in \Pi, \lambda \in \Lambda^*_\varpi \cap \bar{\s}} v^{(P)}_{\varpi,\lambda}$
where $\s$ 
is the open face of $\tcone $ containing $\zeta$. Hence
by Lemma \ref{GJS6.2} 
 $ k_1^{-1} k_2
\in  [K_\zeta(P),K_\zeta(P)] $, and thus $\calf^{(P)}$ induces a 
continuous bijection
$ \tkimpq  \to
\overline{G/U}^{{\rm aff}} \subseteq {E^{(P)}} $.
Since $K$ is compact and so the  map $
(K  \times \liek^{(P)*}_+)/\approx_{(P)} \to
\lieks \times \liek^{(P)*}$
given by 
$(k,\zeta) \mapsto (Ad^*(k)(\zeta),\zeta))$ is proper,
this continuous bijection is a homeomorphism. \hfill $\Box$

\begin{rem}  If ${K^{(P)}}=T$ and $\zeta \in \liek^{(P)*}_+$ then 
$K_\zeta(P)=K_\zeta$, and
so $X_{{\rm KimplT}}$ is the standard
imploded cross-section $\ximp$ of \cite{GJS}. On the other hand if ${K^{(P)}}=K$ then 
$K_\zeta(P)$ is conjugate to $T$ and $[K_\zeta(P),K_\zeta(P)]$ is trivial
for all $\zeta \in \liek^{(P)*}_+$,
so  $X_{{\rm KimplK}} = T^*K$.
\end{rem}

Of course $\overline{G/U}^{\rm aff}$ inherits a $K \times K^{(P)}$-invariant K\"{a}hler
structure as a complex subvariety of $E^{(P)}$. The subvariety  $\overline{G/U}^{\rm aff}$
(which is in general singular) is stratified by the (finitely many) $G$-orbits in
 $\overline{G/U}^{\rm aff}$, and the $K \times K^{(P)}$-invariant K\"{a}hler structure
 on $E^{(P)}$ restricts to a $K \times K^{(P)}$-invariant symplectic structure on 
 each stratum, which
gives $\overline{G/U}^{\rm aff}$
a stratified symplectic structure. Under the homeomorphism 
$ \tkimpq  \to
\overline{G/U}^{{\rm aff}} $ of Theorem \ref{maint} these strata correspond to the
locally closed subsets 
$$ 
\frac{K \times {\rm Ad}^*(K^{(P)})\sigma 
}{\approx^{K^{(P)}}} \cong K^{(P)} \times_{K_\s \cap K^{(P)}} 
\left( \frac{K \times \sigma 
}{\approx^{K^{(P)}}} \right)$$
$$ \cong K^{(P)} \times_{K_\s \cap K^{(P)}} 
\left( \frac{K \times \sigma 
}{
[K_{\s
(P)
},K_{\s
(P)
}]}  \right)$$
of $\tkimpq$ where $\s \in \Sigma 
$ runs over the open faces of $\liets_+$. So the homeomorphism 
$ \tkimpq  \to
\overline{G/U}^{{\rm aff}} $ of Theorem \ref{maint} induces a stratified $K \times K^{(P)}$-invariant
symplectic structure on the universal $K^{(P)}$-imploded cross-section $\tkimpq$. As in
\cite{GJS} the induced symplectic structure on 
 $$  K^{(P)} \times_{K_\s \cap K^{(P)}} 
\left( \frac{K \times \sigma 
}{
[K_{\s
(P) },K_{\s
(P)}]}  \right)$$
can be described directly, and can be expressed in terms of the symplectic
reduction by the action of the subgroup $[K_{\s(P)},K_{\s(P)}]$ of $K$ on a locally closed symplectic
submanifold
of $T^*K$ (cf. \cite{GJS} $\S$2). 

Using this symplectic structure on $\tkimpq$ we obtain the following corollary.

\begin{corollary} \label{lemimp}
Let $K$ act on a symplectic manifold $X$ with moment map $\mu:X \to \lieks$. Then
the symplectic quotient of $\overline{G/U}^{{\rm aff}} \times X = \tkimpq \times X$
by the diagonal action of $K$ can be identified via $\calf^{(P)}$ with $\ximpq$.
\end{corollary}

\begin{rem}
In particular if $X$ is a projective variety with a linear action of 
the complexification $G$ of $K$, then $\ximpq$ can be identified with 
the GIT quotient of $\overline{G/U}^{{\rm aff}} \times X$
by the diagonal action of $G$. \end{rem}

It follows from Corollary \ref{lemimp} that if $(X,\omega)$ is any symplectic manifold on which $K$ acts
with moment map $\mu:X \to \lieks$ then $\ximpq$ inherits a 
stratified $K \times K^{(P)}$-invariant
symplectic structure
$$  \ximpq =  \bigsqcup_{\sigma \in \Sigma} \frac{\mu^{-1}(\sigma)}{\approx^{K^{(P)}}}$$
\begin{equation} =  \mu^{-1}((\liek^{(P)*}_+)^\circ) \sqcup
 \bigsqcup_{\begin{array}{c}\sigma \in \Sigma\\  \sigma \neq (\liets_+)^\circ \end{array}
 } 
 K^{(P)} \times_{K_\s \cap K^{(P)}} 
\left( \frac{\mu^{-1} (\sigma) }{
[K_{\s(P)},K_{\s(P)}]}  \right) \end{equation}
with strata indexed by the set $\Sigma$ of open faces of $\liets_+$,
which are locally
closed symplectic submanifolds of $\ximpq$.
The induced action of ${K^{(P)}}$ on 
$\ximpq$ preserves this symplectic structure and has a moment map
$$\mu_{\ximpq}:\ximpq \to \liek^{(P)*}_+ \subseteq {\liek^{(P)*}}$$
inherited from the restriction of $\mu$ to $\mu^{-1}(\liek^{(P)}_+)$. 

\begin{rem}
In order to identify $\overline{G/U}^{{\rm aff}}$ with $\tkimpq$ we made the
assumption that $K$ is semisimple and simply connected. However the construction of $\ximpq$ makes sense
whenever $K$ is a compact connected Lie group with a Hamiltonian action on the symplectic manifold $X$, and as in \cite{GJS} we can identify $\overline{G/U}^{{\rm aff}}$ with $\tkimpq$ in this more general
situation by expressing $K$ as the quotient of the product of its centre $Z(K)$ and the universal cover
of $[K,K]$ by a finite central subgroup. We then get an identification of $\ximpq$ with the symplectic quotient of 
$\overline{G/U}^{{\rm aff}} \times X$ by $K$ in the general case.
\end{rem}

\subsection{Wonderful compactifications, symplectic cuts and partial desingularisations}

Recently Paradan \cite{Paradan} has introduced a generalisation of the technique of symplectic
cutting (originally due to Lerman \cite{Lerman}) which is valid for a (not necessarily abelian) compact connected group
$K$ and is motivated by the wonderful compactifications of De Concini and Procesi. He defines a
{\em $K$-adapted} polytope in $\liets$ to be a $W$-invariant Delzant polytope $\calp$ in $\liets$
whose vertices are regular elements of the weight lattice $\Lambda^*$. If $\{\l_1,\ldots,\l_N\}$
are the dominant weights lying in the union of all the closed one-dimensional faces of
$\calp$, then there is a $G \times G$-equivariant embedding of $G= K_c$ into
$$\PP( \bigoplus_{i=1}^N V^*_{\l_i} \otimes V_{\l_i})$$
associating to $g \in G$ its representation on $\bigoplus_{i=1}^N V_{\l_i}$. 
The closure $\calx_{(\calp,K)}$ of the image of $G$ in this
projective space is smooth and has moment map
$$\mu^\calp_{K \times K}: \calx_{(\calp,K)} \to \lieks \times \lieks$$
whose image is
$$\mu^\calp_{K \times K}( \calx_{(\calp,K)}) = \{(Ad^*(k_1)\xi,-Ad^*(k_2)\xi) : \xi \in \calp, k_1, k_2 \in K \}.$$
The symplectic cut $X_{(\calp,K)}$ defined by Paradan of a symplectic
manifold $X$ under a Hamiltonian $K$-action with respect to such a $K$-adapted polytope $\calp$ is given
by the symplectic quotient of $ \calx_{(\calp,K)} \times X$ by $K$, so that if $X$ is a complex projective variety
with a linear $K$-action then $X_{(\calp,K)}$ is the GIT quotient
$$X_{(\calp,K)} = (\calx_{(\calp,K)} \times X)/\!/G$$
where $G=K_c$.
Then $X_{(\calp,K)}$ inherits a Hamiltonian $K$-action with moment map
$\mu^{X_{(\calp,K)}}: X_{(\calp,K)} \to \lieks$ whose image is
$$\mu^{X_{(\calp,K)}}(X_{(\calp,K)}) = \mu(X) \cap Ad^*(K)(\calp).$$
Moreoever if $U_\calp = Ad^*(K)(\calp^\circ)$ where $\calp^\circ$ is the
interior of $\calp$ then $(\mu^{X_{(\calp,K)}})^{-1}(U_\calp)$ is an open dense subset of $X_{(\calp,K)}$
which is $K$-equivariantly diffeomorphic to the open subset 
$\mu^{-1}(U_\calp)$ of $X$. This diffeomorphism is a 
quasi-symplectomorphism in the sense that there is a homotopy of
symplectic forms taking the symplectic form on $(\mu^{X_{(\calp,K)}})^{-1}(U_\calp)$
to the pullback of the symplectic form on $\mu^{-1}(U_\calp)$. 

Recall from \cite{GJS} $\S$7 that if $\calp_\e$ is the polyhedral cone
$-(\e \l_0 + \liets_+)$ where $\l_0$ is a generic element of $\mu(X) \cap \liets_+$ and
$0 < \e < \!< 1$, then the imploded cross-section $\ximp = X_{{\rm KimplT}}$ has a partial 
desingularisation
$$\widetilde{X_{{\rm impl}}} = (\ximp)_{(\calp_\e,T)}$$
which is the symplectic reduction of $\calx_{(-\liets_+,T)} \times \ximp $ at $\e \l_0$. 
Similarly, just as in \cite{GJS}, if $P \geq B$ is a parabolic subgroup of
$G=K_c$ with maximal compact subgroup $K^{(P)} = K \cap P$
and unipotent radical $U$, then 
we can construct a $K\times K^{(P)}$-equivariant desingularisation
$\widetilde{\tkimpq}$ for the universal imploded cross-section $\tkimpq \cong \overline{G/U}^{{\rm aff}}$
and a partial desingularisation $\widetilde{\ximpq}$ for $\ximpq$, 
 which can be identified with
the symplectic quotient of $X \times \widetilde{\tkimpq}$ by the induced action of $K$. Moreover $\widetilde{\tkimpq}$ can be identified as a Hamiltonian $K$-manifold with 
\begin{equation} \widetilde{G/U}^{{\rm aff}} = G \times_P (\overline{L^{(P)} \sum_{\varpi \in \Pi} v_{\varpi}^{(P)}})
   = K \times_{K^{(P)}} (\overline{L^{(P)} \sum_{\varpi \in \Pi} v_{\varpi}^{(P)}})    
\end{equation} 
 where $\overline{L^{(P)} \sum_{\varpi \in \Pi} v_{\varpi}^{(P)}}$ is the closure in $E^{(P)}$
 (or equivalently in the
linear subspace $\bigoplus_{\varpi \in \Pi}  (V_{ \varpi}^{K^{(P)}})^* \otimes
V_{\varpi}^{K^{(P)}} $ of $E^{(P)}$)
 of the $L^{(P)}$-orbit (or equivalently the $P$-orbit) of $ \sum_{\varpi \in \Pi} v_{\varpi}^{(P)}$, and
 the restriction 
to $G \times \overline{L^{(P)} \sum_{\varpi \in \Pi} v_{\varpi}^{(P)}}$ of the multiplication map $G \times E^{(P)} \to E^{(P)}$ induces a birational $G$-equivariant
morphism 
$$p_{U}: \widetilde{G/U}^{{\rm aff}} \to \overline{G/U}^{{\rm aff}} = \tkimpq \subseteq E^{(P)}.$$
It follows from Theorem 3.5 of \cite{Paradan} that
$\overline{L^{(P)} \sum_{\varpi \in \Pi} v^{(P)}_{\varpi}})$ is a nonsingular
subvariety of
$$\bigoplus_{\varpi \in \Pi}  (V_{ \varpi}^{K^{(P)}})^* \otimes
V_{\varpi}^{K^{(P)}} \ \subseteq \ \PP(\CC \oplus
\bigoplus_{\varpi \in \Pi}  (V_{ \varpi}^{K^{(P)}})^* \otimes
V_{\varpi}^{K^{(P)}}).$$
If $\lambda_0 \in \mu(X) \cap \liets_+ \cap \liez^{(P)*} $ is generic and $\e > 0$ is sufficiently close to 0, and if $\omega_\e$ is the K\"{a}hler form
on $G/P$ given by regarding $G/P$ as the coadjoint $K$-orbit through $\e\lambda_0$, then $p_{U}^* \omega_{E^{(P)}} + q_P^* \omega_\e$ is a K\"{a}hler form on $\widetilde{G/U}^{{\rm aff}}$ where $q_P: G \times_P E^{(P)} \to G/P$ is the
projection.

The partial desingularisation $\widetilde{\ximpq}$ can alternatively be obtained from
$\ximpq$ via a symplectic cut
following Paradan \cite{Paradan}. Let $W^{(P)}$ be the Weyl group of 
the compact subgroup $K^{(P)}$ of $K$; then we have
an identification 
 \begin{equation} \label{llb} \widetilde{\ximpq} = (\ximpq)_{(\calp_\e,K^{(P)})} \end{equation}
where the cut is with respect to the $K^{(P)}$-action and 
the polyhedral cone $\calp_\e = -(\e \lambda_0 + \tcone)$. If we wish we can cut with respect to 
a suitable $W^{(P)}$-invariant
Delzant polytope  $\calp_\e$  in  this cone
which is large enough that its complement does not meet the
compact subset $\mu(X)$, but then the identification (\ref{llb})
is not quite symplectic according to Paradan's construction; as in Remark 3.1 we
have to distinguish between the flat K\"{a}hler metric on 
$$\bigoplus_{\varpi \in \Pi}  (V_{ \varpi}^{K^{(P)}})^* \otimes
V_{\varpi}^{K^{(P)}}  \subseteq E^{(P)}$$
and the Fubini-Study metric on
$$\bigoplus_{\varpi \in \Pi}  (V_{ \varpi}^{K^{(P)}})^* \otimes
V_{\varpi}^{K^{(P)}} \ \subseteq \ \PP(\CC \oplus
\bigoplus_{\varpi \in \Pi}  (V_{ \varpi}^{K^{(P)}})^* \otimes
V_{\varpi}^{K^{(P)}})\ \subseteq \ \PP(\CC \oplus
E^{(P)}).$$


\section{Non-reductive geometric invariant theory}

The last section discussed a generalisation of symplectic implosion which is closely related to a
GIT-like quotient construction for a linear action of the unipotent radical $U$ of a parabolic subgroup $P$ of a complex reductive group $G$ on a complex variety $X$. This section
will recall from \cite{DK} a version of GIT for non-reductive group actions
 and then relate it to symplectic implosion.

\subsection{Background}

Let $H$ be an affine algebraic 
group, with unipotent radical $U$ (that is, $U$ is the unique maximal normal
unipotent subgroup of $H$), acting linearly on a complex projective variety $X$ with respect to an ample line bundle $L$. If we wish to generalise Mumford's GIT to this non-reductive situation, the first problem to be faced
is that the ring of invariants 
$${\hat{\calo}}_L(X)^H = \bigoplus_{k \geq 0} H^0(X, L^{\otimes k})^H$$
is not necessarily finitely generated as a graded complex algebra,
so that $\Proj({\hat{\calo}}_L(X)^H)$ is not well-defined as a projective variety.
Note, however, that in the case considered in $\S$3 when the unipotent radical $U$ of a parabolic subgroup of a reductive group $G$ acts linearly on $X$ and the linear action extends to
$G$, then the ring of invariants is finitely generated.
Even when $\hat{\calo}_L (X)^H$ is not finitely generated $\Proj({\hat{\calo}}_L(X)^H)$ does make sense as a scheme, and the
inclusion of ${\hat{\calo}}_L(X)^H$ in ${\hat{\calo}}_L(X)$ gives us a rational map of schemes 
$q$ from $ X$ to $ \Proj({\hat{\calo}}_L(X)^H)$, whose
image in   $\Proj({\hat{\calo}}_L(X)^H)$ is constructible (that is, a finite union of
locally closed subschemes). 

We will only consider the case when $H=U$ is unipotent, since $H/U$ is always reductive and classical GIT
allows us to deal with quotients by reductive groups. A more leisurely introduction to non-reductive GIT and details and proofs 
of the results quoted below can be found in \cite{DK}.

\begin{defn} (See \cite{DK}).
Let $I = \bigcup_{m>0} H^0(X,L^{\otimes m})^U$
and for $f \in I$ let $X_f$ be the $U$-invariant affine open subset
of $X$ where $f$ does not vanish, with ${\calo}(X_f)$ its coordinate ring.
The  (finitely generated) {\em semistable set} of $X$ is  $$ X^{ss} = X^{ss,
fg} =  \bigcup_{f \in I^{fg}} X_f$$ where
$I^{fg}$ consists of $f
\in I$ such that ${\calo}(X_f)^U$
 is finitely generated.
 The set of (locally trivial) {\em stable} points is $$  X^s=X^{lts} =
\bigcup_{f \in I^{lts} } X_f$$ where
$I^{lts}$ is the set of $f
\in I$ such that ${\calo}(X_f)^U$
is finitely generated, and $ q: X_f \longrightarrow \Spec({\calo}(X_f)^U)$ is a locally trivial
geometric quotient.
The set of {\em naively semistable}
points of $X$ is the domain of definition 
     $$X^{nss} = \bigcup_{f \in I} X_f$$
of the rational map $q$, and the set of {\em naively stable}
points of $X$ is
     $$X^{ns} = \bigcup_{f \in I^{ns}} X_f$$ where
$I^{ns}$ consists of those $f
\in I$ such that $ {\calo}(X_f)^U$
 is finitely generated, and $  q: X_f \longrightarrow
\Spec({\calo}(X_f)^U)$ is a geometric quotient.


\label{defn:envelopquot}
The {\em enveloped quotient} of
$X^{ss}$ is $q: X^{ss} \rightarrow q(X^{ss})$, where
$q(X^{ss})$ is a dense constructible subset (but not necessarily a subvariety) of the {\em
enveloping quotient}
$$X /\!/ U = \bigcup_{f \in I^{ss,fg}}
\Spec({\calo}(X_f)^U)$$ of $X^{ss}$. 
\end{defn}

\begin{lemma} {\rm (\cite{DK} 4.2.9 and 4.2.10).} \label{patching2} The enveloping quotient
$X/\!/U$ is a quasi-projective variety, and if ${\hat{\calo}}_L(X)^U$ is finitely generated then it is the projective
variety $\Proj({\hat{\calo}}_L(X)^U)$.
\end{lemma}

Let $G$ be a complex reductive group with  $U$
as a closed subgroup, and let $G \times_U X$ denote the quotient of $G \times X$
by the free action of $U$ defined by $h(g,x)=(g h^{-1}, hx)$,
which is a quasi-projective variety by \cite{PopVin} Theorem 4.19. There
is an induced $G$-action on $G \times_U X$ given by left
multiplication of $G$ on itself.
If the action of $U$ on $X$ extends to an action of $G$ there is an isomorphism of
$G$-varieties 
\begin{equation} \label{9Febiso} G \times_U X \cong (G/U) \times X
 \end{equation} given by
$[g,x] \mapsto (gH, gx).$
When $U$ acts linearly on $X$ with respect to a very ample line bundle $L$ inducing an
embedding of $X$ in $\PP^n$, and $G$ is a subgroup of $SL(n+1; \CC)$,
then there is a very ample $G$-linearisation 
(which we will also denote by $L$) on $G \times_U X$
via the embedding
$$ G \times_U X \hookrightarrow G \times_U \mathbb{P}^n  \cong
(G/U) \times \mathbb{P}^n,
$$and using the trivial bundle on the  variety
$G/U$ which is quasi-affine by \cite{Grosshans} Corollary 2.8. For large enough $m$ we can choose a $G$-equivariant embedding of $G/U$ in  $\mathbb{C}^m$ with a linear $G$-action to get a $G$-equivariant
embedding of $G \times_U X$ in
$\mathbb{C}^m \times \mathbb{P}^n \subseteq \mathbb{P}^m \times \mathbb{P}^n
\subseteq \mathbb{P}^{nm+m+n}$
and the $G$-invariants on $ G \times_U X$ are given by
\begin{equation} \label{name}  \bigoplus_{m \geq 0} H^0(G \times_U X, L^{\otimes m})^G \cong
\bigoplus_{m \geq 0} H^0(X, L^{\otimes m})^U = {\hat{\calo}}_L(X)^U.\end{equation}

\begin{defn} \label{defn:separ} (\cite{DK} $\S$5.2).
 A {\em finite separating
set of invariants} for the linear action
of $U$ on $X$ is a collection of invariant sections $\{f_1,
\ldots, f_n \}$ of positive tensor powers of $L$ such that, if $x,y$
are any two points of $X$ then $ f(x) = f(y)$ for all
invariant sections $f$
of $L^{\otimes k}$ and all $k
> 0$ if and only if
$$ f_i(x) =
f_i(y) \hspace{.5in} \forall i = 1, \ldots, n.
$$
If $G$ is any reductive group containing
$U$, a finite separating set $S$ 
of invariant sections of positive tensor powers of $L$ is a {\em
finite fully separating set of invariants} for the linear $U$-action
on $X$ if

(i) for every $x \in X^{s}$ there exists $f \in S$ with associated
$G$-invariant $F$ over $G \times_U X$ (under the isomorphism
(\ref{name})) such that $x \in (G \times_U
X)_{F}$ and $(G \times_U X)_{F}$ is affine; and

(ii) for every $x \in X^{ss}$ there exists $f \in S$ such that $x
\in X_f$ and $S$ is a generating set for ${\calo}(X_f)^U$.

\noindent By \cite{DK} Remark 5.2.3 this
definition is in fact independent of the choice of $G$.

 A
$G$-equivariant projective completion $\overline{G \times_U X}$ of
$G \times_U X$, together with a $G$-linearisation with respect to
a line bundle $L$ which
restricts to the given $U$-linearisation on $X$,
is a
{\em reductive envelope} of the linear $U$-action on $X$ if every
$U$-invariant $f$ in some finite fully separating set of invariants
$S$ for the $U$-action on $X$  extends to a $G$-invariant section of a tensor
power of $L$ over $\overline{G \times_U X}$.
If 
 $L$ is ample on 
$(\overline{G \times_U X})$ 
 it is an {\em ample
reductive envelope}.

\end{defn}

There always exists an ample 
 reductive envelope 
for any linear $U$-action on a projective variety $X$, at least if we replace the line bundle $L$ with a suitable
positive tensor power of itself (see \cite{DK} Proposition
5.2.8).

\begin{defn}\label{defn:s/ssbar} Let $X$ be a
projective variety with a linear $U$-action and a reductive envelope
$\overline{G \times_U X}$. 
Let
 $i: X \hookrightarrow G
\times_U X$ and $j: G \times_U X \hookrightarrow \overline{G
\times_U X}$ be the inclusions, and $\overline{G \times_U X}^{s}$
and $\overline{G \times_U X}^{ss}$ be the stable and semistable sets for
the linear $G$-action on $\overline{G \times_U X}$. 
 Then the set of {\em completely stable points} of $X$ with
respect to the reductive envelope is $$X^{\overline{s}} = (j
\circ i)^{-1}(\overline{G \times_U X}^s)$$ and the set of {\em
completely semistable points} is $$X^{\overline{ss}} = (j \circ
i)^{-1}(\overline{G \times_U X}^{ss}),$$

\end{defn}

\begin{theorem}\label{thm:main} {\rm (\cite{DK} 5.3.1).} 
Let $X$ be a normal projective variety with a linear $U$-action, for $U$ a
connected unipotent group, and let $(\overline{G \times_U X},L)$ be
any ample reductive envelope. Then
there is a
diagram
$$ \begin{array}{ccccccccccccc}
X^{\overline{s}} & \subseteq & X^{s} &
\subseteq & X^{ns} & \subseteq & X^{ss} & \subseteq & X^{\overline{ss}} = X^{nss}\\
\downarrow &  & \downarrow &  & \downarrow & &
\downarrow & & \downarrow \\ X^{\overline{s}}/U & \subseteq &
X^{s}/U & \subseteq & X^{ns}/U & \subseteq & X/\!/U & \subseteq &
\overline{G \times_U X}/\!/G
\end{array}
$$
where all the inclusions are open 
and the first three vertical maps provide quasi-projective geometric quotients of the stable sets $X^{\overline{s}}$,
$X^{s}$ and $X^{ns}$ by the action of $U$. 
 The fourth vertical map 
is the enveloping quotient $q: X^{ss} \rightarrow X /\!/ U$ 
defined in Definition \ref{defn:envelopquot} and
  $X /\!/ U$ is an open
subvariety of the projective variety $\overline{G \times_U X}/\!/G$.
%
\end{theorem}

Note however that, even when  ${\hat{\calo}}_L(X)^U$ is finitely generated so
that  
$$\xu = \Proj({\hat{\calo}}_L(X)^U) = \overline{G \times_U X}/\!/G,$$
the maps $q: X^{ss} \rightarrow X /\!/ U$ and 
$X^{\overline{ss}} \to 
\overline{G \times_U X}/\!/G$ are {\em not} necessarily surjective, and their
images are in general only constructible subsets and not subvarieties.

\subsection{Some examples of reductive envelopes}

Now let us assume that $U = \cplusr$ where $\CC^+$ is the additive group of complex
numbers and $r$ is any positive integer.

\begin{rem}
Each affine algebraic group $H$ over $\CC$ has a unipotent radical $U$, which
is the unique maximal normal unipotent subgroup of $H$ and has a reductive quotient group $R=H/U$ (see e.g. \cite{borel,Springer} for more details).
Given a linear action of $H$ on a projective variety $X$ with respect to an
ample line bundle $L$, we 
can hope to quotient first by the action of $U$, and then
by the induced action of the reductive group $H/U$, provided that
the unipotent quotient (or compactified quotient) is sufficiently canonical to inherit an
induced linear action of $H/U$. 
For example, if the algebra of invariants $\hat{\calo}_L(X)^U$ is finitely
generated then the enveloping quotient $\xu = \Proj(\hat{\calo}_L(X)^U)$
is a projective variety with an induced linear action of $H/U$ on an induced
ample line bundle on $\xu$, and then classical GIT allows us to construct
$X/\!/H =  \Proj(\hat{\calo}_L(X)^H)$ as a GIT quotient $(\xu)/\!/(H/U)$
of $\xu$ by the reductive group $H/U$; even when $\hat{\calo}_L(X)^U$ is not finitely generated, the same is true for $\Proj(\hat{\calo}_L(X)^U_m)$ where $m$ is a sufficiently large positive integer and $\hat{\calo}_L(X)^U_m$ is the subalgebra of $\hat{\calo}_L(X)^U$ generated by invariant sections of $L^{\otimes j}$ for
$1 \leq j \leq m$. 
Moreover the unipotent radical $U$
has canonical sequences of normal subgroups such that each
successive subquotient is isomorphic to $(\CC^+)^r$ for some $r$
(for example by taking the ascending or descending central series of $U$),
so we can hope to quotient successively by unipotent
groups of the form $(\CC^+)^r$, and then finally by the reductive
group $R$. Therefore the case
when $U \cong (\CC^+)^r$ for some $r$ is less special than
it might appear at first sight. \end{rem}

Note that when $U = \cplusr$ we have $\mbox{Aut}(U) \cong \glr$; let
$$\hat{U} = \CC^* \ltimes U$$
be the semidirect product where
 $\CC^*$ is the centre of $\mbox{Aut}(U)$. The centre
of $\hat{U}$ is finite and meets $U$ in the trivial subgroup, so $U$ is isomorphic to a closed subgroup of the reductive
group $G=SL(\CC \oplus \lieu)$
via the inclusion
$$ U \hookrightarrow \hat{U} \to \mbox{Aut}(\hat{U})
\to GL(\mbox{Lie}\hat{U}) = GL(\CC \oplus \lieu)$$
where $\lieu$ is the Lie algebra of $U$ and $\hat{U}$ is identified with its group of inner automorphisms. 
Then $U$ is the unipotent radical of a parabolic
subgroup $P$ of $G = \slrplus$, where $P$ is the stabiliser
of the $r$-dimensional linear subspace $\lieu$ of $\CC \oplus \lieu$,
so we are in the situation of $\S$3.2 above.
The parabolic $P = U   \rtimes  \glr $
 in $G=\slrplus$ has Levi subgroup
$\glr$ embedded in $\slrplus$ as
$$ g \mapsto \left( \begin{array}{ccc}
g & 0 \\ 
0 & \det g^{-1}
\end{array}\right) .$$
 Note that 
$$G/U \cong \{ \alpha \in (\CC^r)^* \otimes \CC^{r+1} |  \ \alpha:\CC^r \to \CC^{r+1} 
\mbox{ is injective }\}$$
with the natural $G$-action $g\alpha = g \circ \alpha $. Since the injective
linear maps from $\CC^r$ to $\CC^{r+1}$ form an open subset in the affine
space $(\CC^r)^* \otimes \CC^{r+1}$ whose complement has codimension two,
we see directly in this case that $U = \cplusr$ is a Grosshans subgroup of $G = \slrplus$
and hence that
$$\calo(G)^U \cong \calo(G/U) \cong \calo((\CC^r)^* \otimes \CC^{r+1})$$
is finitely generated \cite{Grosshans} with 
$$\overline{G/U}^{{\rm aff}} = \Spec \calo(G)^U = (\CC^r)^* \otimes \CC^{r+1}.$$ 

Now suppose that the linear action of $U = \cplusr$
on $X$ extends to a linear action of $G = \slrplus$, giving us an identification of
$G$-spaces
$$G \times_U X \cong (G/U) \times X$$ as at
(\ref{9Febiso}) via $ [g,x] \mapsto (gH, gx)$. Then 
(as in the Borel transfer theorem \cite[Lemma 4.1]{Dolg})
\begin{equation} \hat{\calo}_L(X)^U \cong \hat{\calo}_L(G \times_U X)^G \cong [\calo(G/U) \otimes \hat{\calo}_L(X)]^G \end{equation}
is finitely generated \cite{Grosshans2} and we have a
reductive envelope
$$\overline{G \times_U X} = \prrplus
\times X$$
with 
$$\overline{G \times_U X}/\!/G \cong X/\!/U = \Proj(\hat{\calo}_L(X)^U)$$
where we choose  for our linearisation on $\overline{G \times_U X}$ the line bundle
$$ L^{(N)} = \calo_\prrplus (N) \otimes L$$
with $N>0$ sufficiently large (see \cite{KPEN} $\S$4.1). 
This reductive envelope is ample and so satisfies Theorem \ref{thm:main};
in addition by \cite{KPEN} $\S$4.1(6)
 we have
\begin{equation}
X^{\bar{s}} = X^s \mbox{ and } X^{\bar{ss}} = X^{ss}. \label{sbarssbar} \end{equation}
Thus we have a diagram
$$ \begin{array}{ccccccccccccc}
 X^{s} &
\subseteq & X^{ss} &  & \\
\downarrow &  & \downarrow & &
 \\ 
X^{s}/U & \subseteq & X/\!/U & = &
\overline{G \times_U X}/\!/G
\end{array}
$$
but the enveloping quotient map $q:X^{ss} \to  X/\!/U =
\overline{G \times_U X}/\!/G$ is not necessarily surjective,
so in contrast to the reductive situation we cannot describe
$\xu$ topologically as the quotient of $X^{ss}$ by an equivalence
relation.

In order to describe $\xu$ topologically (and geometrically) it is
useful to consider the linear action of the Levi subgroup $\glr \leq P$
on the closure $\overline{P \times_U X} = \prr \times X$
of $P \times_U X \cong L^{(P)} \times X$ in $\overline{G \times_U X} = \prrplus \times X$.
We have
$$ G \times_U X \cong G \times_P (P \times_U X)$$
where $P/U \cong \glr$ and $G/P \cong \PP^r$ is projective,
so  $G \times_P (\overline{P \times_U X})$
is a projective completion of $G \times_U X$. The induced linearisation
of the action of $G$ on $G \times_P (\overline{P \times_U X})$
is not ample: if we regard  $G \times_P (\overline{P \times_U X})$ as a subvariety in
the obvious way of 
$$G \times_P (\overline{G \times_U X})= G \times_P (\prrplus \times X)
\cong (G/P) \times \prrplus \times X $$
$$\cong \PP^r \times \prrplus \times X$$
then the birational morphism 
$$G \times_P (\overline{P \times_U X}) \to \overline{G \times_U X} \cong \prr \times X$$
given by $[g,y] \mapsto gy$ extends to the projection
$$\PP^r \times \prrplus \times X \to \prrplus \times X$$
and the induced line bundle is the restriction to $G \times_P (\overline{P \times_U X})$ of $\calo_{\prrplus}(N) \otimes L$.
However,  if $\epsilon
\in \QQ \cap (0,\infty)$, the tensor product $\hat{L}_\epsilon = \hat{L}_\epsilon^{(N)}$ 
of this line bundle with the pullback via the morphism
$$G \times_P (\overline{P \times_U X}) \to G/P \cong \PP^r,$$
of the
fractional line bundle $\calo_{\PP^r}(\epsilon)$ provides an ample fractional linearisation
for the action of $G$ on $G \times_P ( \overline{P \times_U X})$ with, 
when $\epsilon$ is sufficiently small, an induced surjective birational
morphism
\begin{equation} \label{xhatu}
\widehat{\xu} =_{df} G \times_P ( \overline{P \times_U X}) /\!/_{\hat{L}_\epsilon} G
\to \overline{G \times_U X} /\!/ G = \xu
\end{equation}
which is an isomorphism over
$$(G \times_U X^{\bar{s}})/G \cong X^{\bar{s}}/U = X^s/U.$$
This line bundle $\hat{L}_\epsilon$ can be thought of as the bundle $G \times_P (\calo_{\prr} (N) \otimes L)$
on $G \times_P (\overline{P \times_U X})$, where now the $P$-action on 
$\calo_{\prr} (N) \otimes L$ is no longer the restriction of the $G$-action on
the line bundle $\calo_{\prrplus} (N) \otimes L$ but has been twisted by $\epsilon$ times
the character of $P$ which restricts to the determinant on $\glr$.

Since $\glr = P/U$ has a central one-parameter subgroup $\CC^*$ we
can modify the linearisation of any linear actions of $P$ and $\glr$ by 
multiplying by $\e$ times the standard character $\det$ of $\glr$ for
any $\e\in \QQ$. By the
Hilbert-Mumford criteria (Proposition \ref{sss} above) we have
\begin{equation} \label{pslss}
\overline{P \times_U X}^{ss,P,\epsilon} \subseteq 
\overline{P \times_U X}^{ss,\glr,\epsilon} \subseteq
\overline{P \times_U X}^{ss,SL(r;\CC)} \end{equation}
where $\overline{P \times_U X}^{ss,\glr,\epsilon}$
and $ \overline{P \times_U X}^{ss,SL(r;\CC)}$ (independent of
$\epsilon$) denote the $\glr$-semistable and $SL(r;\CC)$-semistable sets
of $\overline{P \times_U X} $ after twisting the linearisation
by $\epsilon$ times the character $\det$ of $\glr$; this character
is of course trivial on $SL(r;\CC)$.
It turns out (see \cite{KPEN} $\S$4.1(11))
 that if $\epsilon$ is chosen appropriately (close to $-N/2$ where
$N$ is as in the choice of linearisation above) then
\begin{equation} \overline{P \times_U X}^{ss,\glr,\epsilon} =
(\prr \times X)^{ss,\glr,\epsilon} = \glr \times X \end{equation}
and so quotienting we get
\begin{equation} \overline{P \times_U X}/\!/_{\hat{L}^{(N)}_{-N/2}} \glr \cong X.\end{equation} 
 Therefore
\begin{equation} \label{definecalx}
\calx =_{{\rm def}} \overline{P \times_U X}/\!/ SL(r;\CC)
= (\prr \times X)/\!/ SL(r;\CC)
\end{equation}
 is a projective variety with
a linear action of $\CC^* = \glr/SL(r;\CC)$ which we can twist by $\epsilon$ times the
standard character of $\CC^*$, such that when $\epsilon = -N/2$ we get
\begin{equation} \label{calxquot}
\calx /\!/_{-N/2} \CC^* \ \cong \  X
\end{equation}
while for $\epsilon >0$ sufficiently small  we have a surjection
from an open subset
$(\calx/\!/_{\epsilon} \CC^*)^{\hat{ss}}$ of $\calx /\!/_{\epsilon} \CC^*$ onto
$\widehat{\xu}$, and hence onto $\xu$ (see \cite{KPEN} Proposition 4.6).
More precisely let $(\calx/\!/_{\epsilon} \CC^*)^{\hat{s}}$ be the open
subset 
$\overline{P \times_U X}^{s,P,\epsilon}/\glr$ of
$$\overline{P \times_U X}^{s,\glr,\epsilon}/\glr = 
(\overline{P \times_U X}^{s,\glr,\epsilon}/SL(r;\CC)/\CC^*
= \calx^{s,\epsilon}/\CC^* \subseteq
\calx/\!/_{\epsilon} \CC^*$$
and let $\calx^{\hat{ss},\epsilon} = \pi^{-1}((\calx/\!/_{\epsilon} \CC^*)^{\hat{ss}})$
and $\calx^{\hat{s},\epsilon} = \pi^{-1}((\calx/\!/_{\epsilon} \CC^*)^{\hat{s}})$
where $\pi: \calx^{ss,\epsilon} \to \calx/\!/_{\epsilon} \CC^*$ is the quotient
map, so that
\begin{equation} (\calx/\!/_{\epsilon} \CC^*)^{\hat{s}} = \calx^{\hat{s},\epsilon}/\CC^*. \end{equation}

In this construction we can replace the compactification $\prr$ of $\glr$ by its
{ wonderful compactification} $\widetilde{\prr}$ given by blowing up
$\prr= \{[z:(z_{ij})_{i,j=1}^r]\}$ along the (proper transforms of the) subvarieties 
defined by
$$z=0 \mbox{ and } {\rm rank}(z_{ij}) \leq \ell$$
for $\ell = 0,1,\ldots,r$ and by
$${\rm rank}(z_{ij}) \leq \ell$$
for $\ell = 0,1,\ldots,r-1$ \cite{Kausz}. The action of $SL(r;\CC)$ on
$\widetilde{\prr}$, linearised with respect to a small perturbation
of the pullback of $\calo_{\prr}(1)$, satisfies
$$\widetilde{\prr^{ss}} = \widetilde{\prr^s} \mbox{ and } \widetilde{\prr}/\!/SL(r;\CC) \cong \PP^1.$$
If we take $\overline{P \times_U X}$ to be
$$\widetilde{\prr} \times X$$
instead of $\prr \times X$, and define
\begin{equation} \label{newdefncalx}
\tilde{\calx} =  \widetilde{\prr} \times X/\!/ SL(r;\CC),
\end{equation}
then the properties of $\calx$ given above are satisfied by $\tilde{\calx}$,
and if $X$ is nonsingular then  
$$\widetilde{\xu} =_{{\rm def}} G \times_P (\widetilde{\prr} \times X)/\!/G$$ is a partial desingularisation 
of $\xu$ and a compactification of $X^s/U$.
Indeed it is shown in \cite{KPEN} Proposition 4.6 (combined with 
\cite{KPEN} Remark 4.8) that 
if $\epsilon>0$ is sufficiently small then the natural rational map from
$\tilde{\calx}/\!/_{\epsilon} \CC^*$ to $\widetilde{\xu}$ 
restricts to surjective morphisms
$$(\tilde{\calx}/\!/_{\epsilon} \CC^*)^{\tilde{ss}} \to \widetilde{\xu} \to \xu$$
and
$$(\tilde{\calx}/\!/_{\epsilon} \CC^*)^{\tilde{s}} \to X^s/U.$$
Using the theory of variation of GIT \cite{DolgHu, Ress, Thad}, 
we can relate the quotient $\hat{\calx}/\!/_{\epsilon} \CC^*$ 
 to $\tilde{\calx}/\!/_{N/2} \CC^* \cong X$ via a sequence of flips which occur as walls
are crossed between the linearisations corresponding to $\epsilon$ and to $-N/2$. Thus we have a
diagram
\begin{equation} \begin{array}{ccccccccc}
(\tilde{\calx}/\!/_{\epsilon} \CC^*)^{\tilde{s}} & \subseteq & (\tilde{\calx}/\!/_{\epsilon} \CC^*)^{\tilde{ss}} &
\subseteq & \tilde{\calx}/\!/_{\epsilon} \CC^* & \leftarrow - \rightarrow &
X = \tilde{\calx}/\!/_{-N/2} \CC^*  \\
\downarrow & &  \downarrow & & & \mbox{flips} & \\
X^{{s}}/U & \subseteq &
\widetilde{\xu} & & & & \\
|| & & \downarrow & & & & \\
X^s/U & \subseteq & \xu & & & & 
\end{array}
\end{equation}
where the vertical maps are all {surjective}, 
and
the inclusions are all open.

\begin{rem}

The construction of a reductive envelope described here is only valid if the action
of $U = \cplusr$ on $X$ extends to an action of $G = SL( \CC \oplus \lieu)$ (which is a
rather special situation when the ring of invariants $\hat{\calo}_L(X)^U$ is 
always finitely generated). Moreover at least {\em a priori} this
construction may depend on the choice of
the extension of the $U$-action to a $G$-action, although
$\overline{G \times_U X}/\!/G = \xu = {\rm Proj}(\hat{\calo}_L(X)^U)$ depends
only on the linearisation of the $U$-action on $X$. However it is shown in \cite{KPEN} that we
can associate to a linear $U$-action on $X$ a family of projective varieties $Y_m$ 
(one for every sufficiently large positive integer $m$), each of which contains $X$ and has
an action of $G = SL( \CC \oplus \lieu)$
and a $G$-linearisation on an ample line bundle $L_{Y_m}$, which restricts to the given 
linearisation of the $U$-action on $X$ and is such that every $U$-invariant in a
finite fully separating set of $U$-invariants on $X$ extends to a $U$-invariant on $Y_m$.
Then we can embed $X$ in the $G$-variety
$$\prrplus \times Y_m$$ as $\{\iota \} \times X$ where
$\iota \in (\CC^r)^* \otimes \CC^{r+1} \subseteq \prrplus$
is the standard embedding of $\CC^r$ in $\CC^{r+1}$. The
closure of $GX \cong G \times_U X$ in $\prrplus \times Y_m$ will provide
us with a reductive envelope $\overline{G \times_U X}$
(which is however not necessarily ample), and we can study the closures of the images
of $X^s/U$ in $Y_m/\!/U = \overline{Y_m/\!/U}$ and its partial desingularisation
$\widetilde{Y_m/\!/U}$ constructed as above.
\end{rem}

\subsection{Symplectic implosion for $U=\cplusr \leq \slrplus$ actions}

Let $X$ be a complex projective
variety on which the complexification $G=\slrplus$ of $K=SU(r+1)$
acts linearly with respect to a very ample
line bundle $L$, and 
let $U = \cplusr$ be the unipotent radical of the parabolic $P = \glr U$
as in the previous subsection. As before let $T$ be the maximal torus
of $K$ consisting of the diagonal matrices in $K$, and let $B$ be the
upper triangular Borel subgroup of $G$.
In the notation of $\S$3.2 we have $L^{(P)} = \glr$ and $K^{(P)} = U(r)$.
We can identify the Lie algebra $\liek^{(P)}=\lieu(r)$ of $K^{(P)}$
with the product $[\liek^{(P)},\liek^{(P)}] \oplus \liez^{(P)}$ of the Lie algebras of
its  
semisimple part $Q^{(P)} = [K^{(P)},K^{(P)}] = SU(r)$ and
its centre $Z(K^{(P)}) \cong S^1$. If we identify $\liets$ with
$$\{ \z = (\z_1,\ldots,\z_{r+1}) \in \RR^{r+1} \ \ : \ \ \z_1 + \cdots + \z_{r+1} = 0 \}$$
in the usual way so that
$$\liets_+ = \{  \z = (\z_1,\ldots,\z_{r+1}) \in \liets \ \ : \ \ \z_1 \geq\z_2 \geq \cdots \geq \z_{r+1}  \},$$
then 
\begin{equation} \tcone  = \{  \z = (\z_1,\ldots,\z_{r+1}) \in \liets \ \ : \ \ \z_j \geq \z_{r+1} \mbox{ for }j=1,\ldots,r \} \end{equation}
and
\begin{equation} \liez^{(P)*} = \{  \z = (\z_1,\ldots,\z_{r+1}) \in \liets \ \ : \ \ \z_1 =
\cdots = \z_{r}  \}. \end{equation} 
Moreover $\liek^{(P)*}_+$ can be identified with the set of skew-Hermitian matrices
in ${\mathfrak su}(r+1)^*$ of the form
\begin{equation} \label{bono} 
\z = \left( \begin{array}{cc} \xi & 0 \\ 0 & i\lambda_{r+1} \end{array} \right)
\end{equation}
where $\xi$ is a skew-Hermitian $r \times r$-matrix with all its eigenvalues of
the form $i\lambda$ with $\lambda \in \RR$ and $\lambda \geq \lambda_{r+1}$.
If all the eigenvalues $i\lambda$ of $\xi$ satisfy $\lambda > \lambda_{r+1}$
then 
$K_\zeta(P)$ 
is conjugate to $T$ and $[K_\zeta(P),K_\zeta(P)]$ is trivial.
In general ${\rm Ad}^*(K^{(P)})\z$ contains a matrix of the form
\begin{equation} \label{bonof} 
 \left( \begin{array}{cc} \xi & 0 \\ 0 & i\lambda_{r+1}I_j \end{array} \right)
\end{equation}
for some $j \in \{0,1,\ldots,r\}$, where $\xi$ is a 
skew-Hermitian $(r-j) \times (r-j)$-matrix with all its eigenvalues of
the form $i\lambda$ with $\lambda \in \RR$ and $\lambda > \lambda_{r+1}$,
and $I_j$ is the $j \times j$-identity matrix. Then 
$K_\zeta(P)$ is conjugate in $K^{(P)} = U(r)$ to the product of a torus and the unitary group $U(j)$ embedded in $K=SU(r+1)$ as
$$A \mapsto \left( \begin{array}{ccc} I_{r-j} & 0 & 0 \\ 0 & 
A & 0\\ 0 & 0 &
\det A^{-1} \end{array} \right),$$
and the face $\s$ of $\liek^{(P)*}_+$ to which $\zeta$ belongs is determined
by $j$ and the partition $\pi \in \Pi_{r-j}$ given by the eigenvalues
$i\lambda$  of $\zeta$ with $\lambda > \lambda_{r+1}$.
Thus $[K_\zeta(P),K_\zeta(P)] \cong SU(j)$ and the universal $K^{(P)}$-imploded cross-section is
$$\tkimpq = \bigsqcup_{j=0}^r (K \times \liek^{(P)*}_{+,j,\pi})/\approx^{K^{(P)}}
=  (K \times \liek^{(P)*}_{+})^\circ \ \sqcup
\ \bigsqcup_{j=1}^r \bigsqcup_{\pi \in \Pi_{j}} (K \times \liek^{(P)*}_{+,j})/\approx^{K^{(P)}}$$
\begin{equation}  =  (K \times \liek^{(P)*}_{+})^\circ \ \sqcup
\ \bigsqcup_{j=1}^r \bigsqcup_{\pi=(\pi_1,\ldots,\pi_\ell) \in \Pi_{j}} U(r) \times_{U(\pi_1) \times 
\cdots U(\pi_\ell) \times U(j)} \left(
(K \times \liek^{(P)*}_{+,j,\pi})/SU(j) \right). \end{equation}
Here $\liek^{(P)*}_{+,j}$ consists of all $\zeta \in {\mathfrak su}(r+1)^*$ of the form (\ref{bono}) with 
 $\xi$ a skew-Hermitian $r \times r$-matrix with all its eigenvalues of
the form $i\lambda$ with $\lambda \in \RR$ and $\lambda \geq \lambda_{r+1}$
and exactly $j$ of its eigenvalues equal to $i\lambda_{r+1}$,
and $\liek^{(P)*}_{+,j,\pi}$ consists of all $\zeta \in \liek^{(P)*}_{+,j}$ of the form (\ref{bono}) 
 such that the partition of $r-j$ determined by the eigenvalues of $\z$ of
the form $i\lambda$ with $\lambda > \lambda_{r+1}$ is $\pi$. Moreover if
$(k_1,\zeta_1)$ and $(k_2,\z_2)$ lie in $K \times \liek^{(P)*}_{+,j}$
then $(k_1,\zeta_1) \approx^{K^{(P)}} (k_2,\z_2)$ if and only if 
there is some $\kappa \in K^{(P)}$ such that
$$\z_1 = \z_2 = \kappa 
 \left( \begin{array}{cc} \xi & 0 \\ 0 & i\lambda_{r+1}I_j \end{array} \right)
\kappa^{-1} $$
and 
$\kappa^{-1} k_1 k_2^{-1} \kappa \in [K_\z(P),K_\z(P)] 
\cong SU(j)$. Thus $\tkimpq$ is isomorphic to $\overline{G/U}^{{\rm aff}}
= (\CC^r)^* \otimes \CC^{r+1}$ via
$$(k,\z) \mapsto   k \circ \calf(\zeta)$$
where if $\z$ is at (\ref{bono}) then $\calf(\z):\CC^r \to \CC^r \subseteq\CC^{r+1}$ is the linear map represented by the unique $r \times r$-Hermitian positive definite matrix $\a$ satisfying
$i\a^* \a = \xi - i \lambda_{r+1}I_r$.

Let $\omega$ be a $K$-invariant K\"{a}hler form on
$X$, given in some choice of coordinates by the Fubini-Study
form on the projective space into which the very ample line bundle $L$
embeds $X$. Then we know that 
$$\hat{\calo}_L(X)^U \cong (\hat{\calo}_L(X) \otimes \calo(G)^U)^G$$
is finitely generated, and the associated projective variety
$$X/\!/U = {\rm Proj}(\hat{\calo}_L(X)^U)$$
is isomorphic to the GIT quotient $ (\overline{G/U}^{{\rm aff}} \times X)/\!/G$,
which as in $\S$3.2 can be identified with a symplectic quotient of 
$\overline{G/U}^{{\rm aff}} \times X$ by $K$, and thus with the
 { ${K^{(P)}}$-imploded cross-section}
$$ \ximpq = \mu^{-1}(\liek^{(P)*}_+)/\approx_{K^{(P)}} $$
of $X$, where  $x \approx_{K^{(P)}} y$ if and only if 
$\mu(x) = \mu(y) = \zeta \in \liek^{(P)*}_+$
and
$x = \kappa y$ for some $\kappa  \in [K_\zeta(P),K_\zeta(P)]$. 
Equivalently
$$  \ximpq =  \mu^{-1}((\liek^{(P)*}_{+})^\circ) \ \sqcup
\ \bigsqcup_{j=1}^r \mu^{-1}(\liek^{(P)*}_{+,j})/\approx^{K^{(P)}}$$
\begin{equation} =  \mu^{-1}((\liek^{(P)*}_{+})^\circ) \ \sqcup
\ \bigsqcup_{j=1}^r 
\bigsqcup_{\pi=(\pi_1,\ldots,\pi_\ell) \in \Pi_{j}} U(r) \times_{U(\pi_1) \times 
\cdots U(\pi_\ell) \times U(j)} \left(
\mu^{-1}(\liek^{(P)*}_{+,j,\pi} \cap \liets_+)/SU(j) \right) \end{equation}
since $[K_\z(P),K_\z(P)]  \cong SU(j)$
if $\z \in \liek^{(P)*}_{+,j}$.

The desingularisation $\widetilde{\tkimpq}$ of $\tkimpq$ is given by 
\begin{equation} \label{ooo} \widetilde{\tkimpq} = (K \times \liek^{(P)*,\e}_+)/\approx^{K^{(P)}}_\e \end{equation}
where $\liek^{(P)*,\e}_+ = {\rm Ad}^*(K^{(P)})(\e \lambda_0 + \tcone)$ for 
 $0 < \e <\!< 1$ and 
 $\lambda_0 = \diag(1,1, \cdots, 1, -r) \in \tcone \cap \liez^{(P)*}$,
and if 
$(k_1,\zeta_1)$ and $(k_2,\z_2)$ lie in $K \times \liek^{(P)*,\e}_{+,j}$
then $(k_1,\zeta_1) \approx^{K^{(P)}}_\e (k_2,\z_2)$ if and only if 
there is some $\kappa \in K^{(P)} \cong U(r)$ such that
$$\z_1 = \z_2 = \kappa 
 \left( \begin{array}{cc} \xi & 0 \\ 0 & i\lambda_{r+1}I_j \end{array} \right)
\kappa^{-1} $$
and 
$\kappa^{-1} k_1 k_2^{-1} \kappa$ lies in the maximal torus $T_j$ of
$ [K_\z(P),K_\z(P)] \cong SU(j)$ which is its intersection with $T$. The partial desingularisation
$\widetilde{\ximpq}$ of $\ximpq$ is the symplectic quotient
of $\widetilde{\tkimpq} \times X$ by the diagonal action of $K$; as a stratified symplectic space,
it is given by
$$  \widetilde{\ximpq} =  \mu^{-1}((\liek^{(P)*,\e}_{+})^\circ) \ \sqcup
\ \bigsqcup_{j=1}^r 
\bigsqcup_{\pi=(\pi_1,\ldots,\pi_\ell) \in \Pi_{j}} U(r) \times_{U(\pi_1) \times 
\cdots U(\pi_\ell) \times U(j)} \left(
\mu^{-1}(\e \lambda_0 + \liek^{(P)*}_{+,j,\pi} \cap \liets_+)/T_j \right)$$
and it can also be identified with the partial desingularisation $\widetilde{\xu}$
described in $\S$4.2.

\begin{ex} \label{lastex}
Let $U =\CC^+$ act linearly on a projective space $\PP^n$, and suppose that 
coordinates have been chosen so that the natural generator of ${\rm Lie}(\CC^+) = \CC$ has Jordan normal 
form with blocks
of sizes $k_1 + 1, \ldots, k_s + 1$ where $\sum_{j=1}^s (k_j +1) = n+1.$
The $\CC^+$ action extends to an action of $G= SL(2;\CC)$ by identifying
$ \CC^+$ with the group of upper triangular matrices $$ \{ \left( \begin{array}{cccc}
1& a\\
0&1
\end{array} \right) : \ a \in \CC \} \leq SL(2;\CC)$$
and
$\CC^{n+1}$ with $ \bigoplus_{j=1}^s {\rm Sym}^{k_j}(\CC^2)$
where ${\rm Sym}^k(\CC^2)$ is the $k$th symmetric power of the
standard representation $\CC^2$ of $G = SL(2;\CC)$. We have
$$G/\CC^+ \cong \CC^2 \setminus \{ 0 \} \subseteq \CC^2 \subseteq
\PP^2 = \overline{G/\CC^+}$$ and thus 
$\PP^n /\!/ \CC^+$ is the GIT quotient ${\rm Proj } (\CC[x_0, \ldots, x_n]^{\CC^+}) \cong (\PP^2 \times \PP^n)/\!/G$
with respect to the linearisation $\calo_{\PP^2}(N) \otimes \calo_{\PP^n}(1)$
on $\PP^2 \times \PP^n$ for $N$ a sufficiently large positive integer.
Since 
$(\PP^2)^{ss,G} = \CC^2$ and $N$ is large we have
$$(\PP^2 \times \PP^n)^{ss,G} \subseteq \CC^2 \times \PP^n
= (G \times_{\CC^+} \PP^n) \ \sqcup \  (\{ 0\} \times \PP^n)$$
and if semistability implies stability then
$$\PP^n /\!/ \CC^+ = (\PP^n)^{s,U}/\CC^+ \ \sqcup \ (\{0\} \times \PP^n)/\!/SL(2;\CC).$$
In this example the parabolic subgroup $P$ of $G = SL(2;\CC)$ is its standard (upper triangular) Borel subgroup
with $\overline{B/\CC^+} = \overline{\CC^*} = \PP^1$ and
$$\overline{B \times_{\CC^+} \PP^n} = \PP^1 \times \PP^n,$$
while $G \times_B \overline{B/\CC^+} = G \times_B \PP^1$ is the blow-up
of $\PP^2$ at the origin $0 \in \CC^2 \subseteq \PP^2$. Similarly
$G \times_B (\overline{B \times_{\CC^+} \PP^n})$ is the blow-up 
of $\overline{G \times_{\CC^+} \PP^n} \cong \PP^2 \times \PP^n$
along $\{ 0 \} \times \PP^n$, and its quotient $\widetilde{\xu}$ is the blow-up of
$\PP^n /\!/\CC^+$ along its \lq boundary' 
$$\PP^n /\!/SL(2;\CC) \cong (\{0\} \times \PP^n)/\!/SL(2;\CC) \subseteq (\PP^2 \times \PP^n)/\!/SL(2;\CC)
= \PP^n/\!/\CC^+.$$
From the point of view of symplectic geometry we have 
$$\PP^n/\!/\CC^+ \cong (\PP^n)_{{\rm impl}} 
= \mu^{-1}((\liets_+)^\circ ) \sqcup \frac{\mu^{-1}(0)}{SU(2)}
=\mu^{-1}(0,\infty) \sqcup \frac{\mu^{-1}(0)}{SU(2)}$$
where $\liets_+$ is identified with $(0,\infty)$ in the usual way, and 
$$\widetilde{\PP^n/\!/\CC^+} \cong
 \widetilde{(\PP^n)_{{\rm impl}}} 
 = \mu^{-1}(\e,\infty) \sqcup \frac{\mu^{-1}(\e)}{S^1}$$
for $0<\e <\!< 1$.
\end{ex}

\end{document}